\documentclass[10pt,a4paper]{amsart} 
\usepackage[top=15mm,bottom=15mm,left=15mm,right=15mm]{geometry} 
\usepackage{graphicx}
\usepackage{pgfplots}
\pgfplotsset{compat=1.18} 
\usepackage{amsmath,amssymb,amsthm}
\usepackage{mathtools}
\mathtoolsset{showonlyrefs=true}
\usepackage{bm}
\usepackage{comment}
\usepackage{float}
\usepackage{ascmac} 
\usepackage{color} 
\usepackage{cite}
\usepackage{multirow}
\usepackage{mathrsfs} 
\usepackage{lineno}
\usepackage{empheq} 
\usepackage{enumitem} 

\newtheorem{theorem}{Theorem}[section] 
\newtheorem{lemma}[theorem]{Lemma} 
\newtheorem{proposition}[theorem]{Proposition}

\numberwithin{equation}{section} 

\newcommand{\norm}[2]{\left\| #2 \right\|_{#1}} 
\newcommand{\snorm}[2]{\left| #2 \right|_{#1}} 
\newcommand{\abs}[1]{\left| #1 \right|} 

\title[Conservative semi-Lagrangian scheme with higher-order interpolation]{Error analysis of a conservative semi-Lagrangian scheme with higher-order interpolation for the one-dimensional advection--diffusion equation} 
\author{Haruki Takemura} 
\email{36takemura@gmail.com} 
\date{} 


\begin{document} 
\begin{abstract}
  We present error estimates for a semi-Lagrangian scheme for the one-dimensional advection--diffusion equation in conservation form. 
  We consider a fully semi-Lagrangian scheme, in which the semi-Lagrangian technique is applied not only to advection but also to diffusion. 
  As with other fully semi-Lagrangian schemes, it is derived from a representation of the exact solution involving characteristic curves. 
  The scheme possesses a mass-conservation property analogous to that of the original partial differential equation. 
  To achieve mass conservation, we use an interpolation operator of degree $ 2 s $ ($ s \ge 1 $) that preserves the integral of the interpolated function on each grid cell. 
  We show that this interpolation operator is stable with respect to a modified $ H^s $-norm, which we call the weighted $ H^s $-norm, that depends on the mesh size $ h $ and the time step size $ \Delta t $. 
  We use this fact to establish the stability of the entire fully semi-Lagrangian scheme. 
  We also obtain an $ L^2 $-error estimate of order $ O (\Delta t + h^{2 s + 1} / \Delta t) $ and an $ H^s $-error estimate of order $ O (\Delta t^{3 / 2} / h^s + h^{s + 1} / \Delta t^{1 / 2}) $. 
\end{abstract} 
\maketitle

\section{Introduction} \label{s:1} 
The semi-Lagrangian (SL) technique is a numerical scheme based on the method of characteristics and is widely used to solve problems involving advection. 
By tracing the characteristic curves backward in time at each time step, SL schemes are not restricted by the Courant--Friedrichs--Lewy condition. 
Moreover, the SL technique can also be used for diffusion, not only for advection. 
The schemes obtained in this way are called fully SL (FSL) schemes \cite{18BoFeRo,25Fe}, which we describe in detail below. 

Let $ \mathbb{T} = \mathbb{R} / \mathbb{Z} $, $ T > 0 $, and $ Q = \mathbb{T} \times (0, T) $. 
We consider the following one-dimensional advection--diffusion equation with initial condition 
\begin{empheq}[left=\empheqlbrace]{alignat=2} 
  & \partial_t u (x, t) + \partial_x (a (x, t) \, u (x, t)) - \frac{1}{2} \sigma^2 \partial_{xx} u (x, t) = g (x, t) \quad & & \text{for} \ (x, t) \in Q, \label{5} \\ 
  & u (x, 0) = u_0 (x) & & \text{for} \ x \in \mathbb{T}, \label{6} 
\end{empheq} 
where $ \sigma $ is a positive constant. 
We are concerned with FSL schemes for \eqref{5}--\eqref{6}. 
Let $ N_x, N_t $ be positive integers, $ 0 = x_0 < x_1 < \cdots < x_{N_x} = 1 $, $ E_j = (x_{j - 1}, x_j) $, $ h = \max_{1 \le j \le N_x} |x_j - x_{j - 1}| $, $ \Delta t = T / N_t $, and $ t_n = n \Delta t $. 
We define the space of piecewise linear functions by 
\begin{align}
  V_h^{\mathrm{lin}} = \left\{ v \in C (\mathbb{T}) \mathrel{}\middle|\mathrel{} v|_{E_j} \in \mathbb{P}_1 (E_j) \right\}, 
\end{align} 
where $ \mathbb{P}_k (E_j) $ denotes the set of all polynomials of degree at most $ k $ on $ E_j $. 
We define the linear interpolation operator $ I_h^{\mathrm{lin}}: C (\mathbb{T}) \to V_h^{\mathrm{lin}} $ for continuous functions on $ \mathbb{T} $ by 
\begin{align} 
  (I_h^{\mathrm{lin}} v) (x_j) = v (x_j) \quad \text{for} \ j = 1, \ldots, N_x. 
\end{align} 
Before constructing a conservative scheme, we present the FSL scheme with linear interpolation for \eqref{5}--\eqref{6}. This scheme does not, in general, conserve mass and can be written as follows: 
\begin{empheq}[left=\empheqlbrace]{alignat=2} 
  & u_h^0 = I_h^{\mathrm{lin}} u_0, \label{7-2} \\ 
  & u_h^{n + 1} = I_h^{\mathrm{lin}} \left(\frac{1}{2} (u_h^n \circ X_{\Delta t, +}^{n + 1}) \, Y_{\Delta t}^{n + 1} + \frac{1}{2} (u_h^n \circ X_{\Delta t, -}^{n + 1}) \, Y_{\Delta t}^{n + 1} + Z_{\Delta t}^{n + 1}\right), \label{8-2} 
\end{empheq} 
where 
\begin{align} 
  X_{\Delta t, \pm}^n (x) & = x - a (x, t_n) \Delta t \pm \sigma \sqrt{\Delta t}, \label{9-2} \\ 
  Y_{\Delta t}^n (x) & = 1 - \partial_x a (x, t_n) \, \Delta t, \label{11} \\ 
  Z_{\Delta t}^n (x) & = g (x, t_n) \Delta t. \label{12} 
\end{align} 
To improve accuracy, higher-order interpolation may be employed instead of linear interpolation. 
For instance, an SL scheme using cubic spline interpolation was studied in \cite{76Pu}. 
Furthermore, the SL scheme with the cubic Hermite interpolation is called the cubic interpolated pseudo-particle (CIP) scheme \cite{85TaNiYa}. 

We discuss the mass-conservation property of \eqref{5}--\eqref{6}. 
Let $ u^n (x) = u (x, t_n) $. 
Suppose that $ g \equiv 0 $. 
Integrating \eqref{5} over $ \mathbb{T} \times (0, t_n) $, we have 
\begin{align} 
  \int_\mathbb{T} u^n (x) \, dx = \int_\mathbb{T} u_0 (x) \, dx. \label{2-8-2} 
\end{align} 
However, the scheme \eqref{7-2}--\eqref{8-2} does not conserve the total mass. 
Namely, $\int_\mathbb{T} u_h^n (x) \, dx $ is not equal to $ \int_\mathbb{T} u_0(x) \, dx$ in general. 
The same holds for schemes employing cubic spline or Hermite interpolation. 
Without mass conservation, the scheme may become less accurate, particularly for computations over a long time interval. 
To overcome this issue, several conservative SL schemes have been proposed. 
For instance, a quasi-conservative SL scheme and a finite volume SL scheme were developed in \cite{93Pr} and \cite{95ScSe}, respectively. 
Furthermore, conservative CIP schemes have been proposed for several problems, including Burgers' equation \cite{00TaNaYa, 01YaTaNaXi} and the shallow water equations \cite{09ToOgYa}. 
In the conservative CIP schemes, instead of cubic Hermite interpolation, piecewise quartic polynomial interpolation is employed to achieve the mass-conservation property and a high-order convergence rate simultaneously. 

In this paper, we consider a scheme combining the FSL scheme and the conservative CIP scheme. 
Following \cite{00TaNaYa}, we define an interpolation operator that conserves the integral of the function in each cell $ E_j $ ($ j = 1, \ldots, N_x $). 
For a positive integer $ s $, let 
\begin{align} 
  V_h^{s - 1, 2 s} = \left\{ v \in C^{s - 1} (\mathbb{T}) \mathrel{}\middle|\mathrel{} v|_{E_j} \in \mathbb{P}_{2s} (E_j) \right\}, 
\end{align} 
and let the interpolation operator $ \Pi_h^{2s}: C^{s - 1} (\mathbb{T}) \to V_h^{s - 1, 2 s} $ be defined by 
\begin{align} 
  \partial_x^k (\Pi_h^{2s} v)(x_j) = \partial_x^k v(x_j) \quad (k = 0,\ldots, s - 1, \ j = 1, \ldots, N_x) \label{6-5} 
\end{align} 
and 
\begin{align} 
  \int_{E_j} (\Pi_h^{2s} v)(x) \, dx = \int_{E_j} v(x) \, dx \quad  (j = 1, \ldots, N_x) 
  \label{6-1} 
\end{align} 
for any $ v \in C^{s - 1} (\mathbb{T}) $. 
We note that the interpolation used in \cite{00TaNaYa} is identical to $ \Pi_h^4 $. 
We replace the linear interpolation operator $ I_h^{\mathrm{lin}} $ in \eqref{7-2}--\eqref{8-2} with $ \Pi_h^{2s} $ to obtain the following conservative FSL scheme, which is the main focus of this paper: 
\begin{empheq}[left=\empheqlbrace] {alignat=2} 
  & u_h^0 = \Pi_h^{2s} u_0, \label{3} \\ 
  & u_h^{n + 1} = \Pi_h^{2s} \left(\frac{1}{2} (u_h^n \circ X_{\Delta t, +}^{n + 1}) \, Y_{\Delta t}^{n + 1} + \frac{1}{2} (u_h^n \circ X_{\Delta t, -}^{n + 1}) \, Y_{\Delta t}^{n + 1} + Z_{\Delta t}^{n + 1}\right). \label{4} 
\end{empheq} 
The scheme \eqref{3}--\eqref{4} possesses the mass-conservation property as stated in the following proposition. 

\begin{proposition} \label{p:1} 
  Let $ u_h^n $ be the numerical solution of the scheme \eqref{3}--\eqref{4}. 
  Then, the following holds for $ n = 1, \ldots, N_t $: 
  \begin{align} 
    \int_{\mathbb{T}} u_h^{n} (x) \, dx 
    = \int_\mathbb{T} u_0 (x) \, dx + \Delta t \sum_{k = 1}^{n} \int_{\mathbb{T}} g^k (x) \, dx, \label{2-8} 
  \end{align} 
  where $ g^k (x) = g (x, t_k) $. 
\end{proposition} 
The proof of Proposition~\ref{p:1} is given in Section~\ref{s:2}. 
We note that, when $ g \equiv 0 $, \eqref{2-8} corresponds to \eqref{2-8-2}. 

We briefly recall error estimates of the SL and FSL schemes. 
For the error estimates of FSL schemes with linear interpolation, we refer to \cite{02Mi}, which considers an advection-diffusion equation in non-conservative form with a velocity $ a $ depending on $ u $. 
An $ L^\infty $-error estimate of the finite volume SL scheme proposed in \cite{95ScSe} is given in \cite{01PhWi} for the one-dimensional case. 
Error estimates of the CIP scheme for the one-dimensional linear advection equation are given in \cite{26KaTa}. 

We state the main result of this paper, error estimates of the scheme \eqref{3}--\eqref{4}, as follows. 

\begin{theorem} \label{t:1} 
  For a positive integer $ s $, let $ r = \max \{2s + 1, s + 4\} $, and let $ u \in W^{2, \infty} (0, T; W^{r, \infty} (\mathbb{T})) $ be the unique solution to \eqref{5}--\eqref{6} with 
  \begin{align} 
    a \in W^{1,\infty}(0,T;W^{s+1,\infty}(\mathbb T)), \quad 
    g \in W^{1,\infty}(0,T;W^{s,\infty}(\mathbb T)). 
  \end{align} 
  Suppose that the spatial mesh size $ h $ and the time step size $ \Delta t $ satisfy 
  \begin{gather} 
    0 < \Delta t \le 1, \label{l2-14-1} \\ 
    h^{2s} \le \kappa \Delta t, \label{l2-14} \\ 
    \Delta t \norm{L^{\infty} (Q)}{\partial_x a} \le \frac{1}{2} \label{l2-14-3} 
  \end{gather} 
  for a positive constant $ \kappa $. 
  Let $ \{u_h^n\}_{n = 0}^{N_t} $ be the numerical solution of the scheme \eqref{3}--\eqref{4}. 
  Then, there exists a positive constant $ C $ depending on $ T $, $ s $, 
  $ \norm{W^{2, \infty} (0, T; W^{r, \infty} (\mathbb{T}))}{u} $, 
  $ \norm{W^{1, \infty} (0, T; W^{s + 1, \infty} (\mathbb{T}))}{a} $, 
  $ \norm{W^{1, \infty} (0, T; W^{s, \infty} (\mathbb{T}))}{g} $, 
  $ \sigma $, and $ \kappa $ such that for any $ n = 0, \ldots, N_t $ 
  \begin{align} 
    \norm{L^2 (\mathbb{T})}{u^n - u_h^n} & \leq C \left(\Delta t + \frac{h^{2s + 1}}{\Delta t}\right), \label{t1-1} \\ 
    \norm{H^s (\mathbb{T})}{u^n - u_h^n} & \leq C \left(\frac{\Delta t^{3 / 2}}{h^s} + \frac{h^{s + 1}}{\Delta t^{1 / 2}}\right), \label{t1-2} 
  \end{align} 
  where $ u^n (x) = u (x, t_n) $. 
\end{theorem} 

A key step in this error analysis is to introduce the weighted $H^s$-norm defined by \eqref{3-4}. 
First, we show that the interpolation operator $\Pi_h^{2s}$ is stable with respect to this weighted $H^s$-norm. 
Using the stability of $ \Pi_h^{2s} $, we establish the stability of the scheme in the same norm. 
Finally, by estimating the truncation error, we obtain the error estimates stated in Theorem~\ref{t:1}. 

This paper is organized as follows. 
In Section~\ref{s:2}, we derive the FSL scheme \eqref{3}--\eqref{4} from a stochastic representation of the solution to \eqref{5}--\eqref{6}, and give the proof of Proposition~\ref{p:1}. 
In Section~\ref{s:3}, we give the proof of Theorem~\ref{t:1}. 
In Section~\ref{s:e}, we present some numerical results. 
In Section~\ref{s:a}, we prove an auxiliary lemma used in Section~\ref{s:3}. 

\section{Derivation of the conservative scheme} \label{s:2} 
We denote the norm of the space $ L^2 (\mathbb{T}) $ by $ \norm{0, 2}{\, \cdot \,} $. 
Let $ s $ be a positive integer. 
We also use the Sobolev seminorm $ \snorm{s, 2}{v} = \norm{0, 2}{\partial_x^s v} $. 
Furthermore, for $ \Delta = (h, \Delta t) $, we introduce a norm 
\begin{gather} 
  \norm{s, 2, \Delta}{v} = \left(\norm{0, 2}{v}^2 + \frac{h^{2s}}{\Delta t} \snorm{s, 2}{v}^2 \right)^{1/2}, \label{3-4} 
\end{gather} 
which we refer to as the \emph{weighted $ H^s $-norm} in this paper. 
We also use a generic positive constant $ C $, independent of $ h $ and $ \Delta t $, which may change from line to line. 

\subsection{Derivation of FSL schemes} 
We derive the FSL scheme \eqref{7-2}--\eqref{8-2} and the conservative FSL scheme \eqref{3}--\eqref{4}. 
By the Feynman--Kac formula, the solution to \eqref{5}--\eqref{6} can be expressed as follows (see Chapter~2 of \cite{85Fr}): 
\begin{align} 
  u (x, t) & = \mathbb{E} \bigl[u_0 (X (t; x, t)) \, Y (t; x, t) + Z (t; x, t)\bigr], 
\end{align} 
where $ X (s; x, t) $, $ Y (s; x, t) $, and $ Z (s; x, t) $ are the solution to the stochastic differential equations 
\begin{empheq}[left=\empheqlbrace]{alignat=2} 
  & dX (s; x, t) = - a (X (s; x, t), t - s) \, ds + \sigma \, d w (s) \quad & & \text{for} \ s \in (0, t), \\ 
  & dY (s; x, t) = - \partial_x a (X (s; x, t), t - s) \, Y (s; x, t) \, ds \quad & & \text{for} \ s \in (0, t), \\ 
  & dZ (s; x, t) = g (X (s; x, t), t - s) \, Y (s; x, t) \, ds \quad & & \text{for} \ s \in (0, t), \\ 
  & X (0; x, t) = x, \ Y (0; x, t) = 1, \ Z (0; x, t) = 0. 
\end{empheq} 
Here $ w (s) $ is the Wiener process. 
We can express $ u^{n + 1} $ in terms of $ u^{n} $ as 
\begin{align} 
  u^{n + 1} (x) & = \mathbb{E} \bigl[u^n (X (\Delta t; x, t_{n + 1})) \, Y (\Delta t; x, t_{n + 1}) + Z (\Delta t; x, t_{n + 1})\bigr]. \label{1} 
\end{align} 
Let $ \xi $ be a random variable taking values $ 1 $ and $ - 1 $ with probability $ 1 / 2 $ each. 
We approximate $ X (\Delta t; \,\cdot\,, t_n) $ by a stochastic mapping $ X_{\Delta t}^n $ from $ \mathbb{T} $ to $ \mathbb{T} $ defined by 
\begin{align} 
  X_{\Delta t}^n (x) = x - a (x, t_n) \Delta t + \xi \sigma \sqrt{\Delta t}. 
\end{align} 
In addition, as approximations of $ Y (\Delta t; x, t_n) $ and $ Z (\Delta t; x, t_n) $, we introduce $ Y_{\Delta t}^n (x) $ and $ Z_{\Delta t}^n (x) $ defined in \eqref{11} and \eqref{12}, respectively. 
We discretize \eqref{1} in time to obtain the semi-discrete scheme 
\begin{align} 
  \tilde{u}^{n + 1} & = \mathbb{E} \bigl[(\tilde{u}^n \circ X_{\Delta t}^{n + 1}) \, Y_{\Delta t}^{n + 1} + Z_{\Delta t}^{n + 1} \bigr] \\ 
  & = \frac{1}{2} (\tilde{u}^n \circ X_{\Delta t, +}^{n + 1}) \, Y_{\Delta t}^{n + 1} + \frac{1}{2} (\tilde{u}^n \circ X_{\Delta t, -}^{n + 1}) \, Y_{\Delta t}^{n + 1} + Z_{\Delta t}^{n + 1}, \label{14} 
\end{align} 
where $ X_{\Delta t, +}^n $ and $ X_{\Delta t, -}^n $ are defined by \eqref{9-2}. 

Applying the linear interpolation operator to the right-hand side of \eqref{14}, we obtain the FSL scheme \eqref{7-2}--\eqref{8-2}. 
Applying $ \Pi_h^{2s} $, rather than linear interpolation, we obtain the conservative FSL scheme \eqref{3}--\eqref{4}. 

\subsection{Mass-conservation property of the scheme} 
We now give the proof of Proposition~\ref{p:1}. 

\begin{proof}[Proof of Proposition~\ref{p:1}] 
  Summing \eqref{6-1} from $ j = 1 $ to $ N_x $, we have 
  \begin{align} 
    \int_{\mathbb{T}} (\Pi_h^{2s} v) (x) \, dx = \int_{\mathbb{T}} v (x) \, dx \quad \text{for} \ v \in C^{s - 1}(\mathbb{T}). \label{2-4} 
  \end{align} 
  Integrating \eqref{4} over $ \mathbb{T} $ and using \eqref{2-4}, we obtain 
  \begin{align} 
    \int_{\mathbb{T}} u_h^{n + 1} (x) \, dx & = \int_{\mathbb{T}} \Pi_h^{2s} \left(\frac{1}{2} (u_h^n \circ X_{\Delta t, +}^{n + 1}) \, Y_{\Delta t}^{n + 1} + \frac{1}{2} (u_h^n \circ X_{\Delta t, -}^{n + 1}) \, Y_{\Delta t}^{n + 1} + Z_{\Delta t}^{n + 1}\right) \, dx \\ 
    & = \int_{\mathbb{T}} \left(\frac{1}{2} (u_h^n \circ X_{\Delta t, +}^{n + 1}) \, Y_{\Delta t}^{n + 1} + \frac{1}{2} (u_h^n \circ X_{\Delta t, -}^{n + 1}) \, Y_{\Delta t}^{n + 1} + Z_{\Delta t}^{n + 1}\right) \, dx. \label{2-9} 
  \end{align} 
  Let $ d \in \mathbb{R} $. 
  By the change of variables $ x \mapsto y = x - a (x, t_{n + 1}) \Delta t + d $, we have for any $ v \in C^{s - 1}(\mathbb{T}) $ 
  \begin{align} 
    & \int_{\mathbb{T}} v (x - a (x, t_{n + 1}) \Delta t + d) \cdot Y_{\Delta t}^{n + 1} (x) \, dx \\ 
    & = \int_{\mathbb{T}} v (x - a (x, t_{n + 1}) \Delta t + d) \cdot \partial_x (x - a (x, t_{n + 1}) \Delta t + d) \, dx \\ 
    & = \int_{\mathbb{T}} v(y) \, dy. 
  \end{align} 
  Substituting $ v = u_h^{n} $ and $ d = \pm \sigma \sqrt{\Delta t} $, we have 
  \begin{align} 
    \int_{\mathbb{T}} u_h^n (X_{\Delta t, +}^{n + 1} (x)) \, Y_{\Delta t}^{n + 1} (x) \, dx & = \int_{\mathbb{T}} u_h^n (x) \, dx, \label{2-1} \\ 
    \int_{\mathbb{T}} u_h^n (X_{\Delta t, -}^{n + 1} (x)) \, Y_{\Delta t}^{n + 1} (x) \, dx & = \int_{\mathbb{T}} u_h^n (x) \, dx. \label{2-2} 
  \end{align} 
  Applying \eqref{2-1} and \eqref{2-2} to \eqref{2-9}, we have 
  \begin{align} 
    \int_{\mathbb{T}} u_h^{n + 1} (x) \, dx & = \int_{\mathbb{T}} u_h^{n} (x) \, dx + \int_{\mathbb{T}} Z_{\Delta t}^{n + 1} (x) \, dx. \label{2-5} 
  \end{align}  
  Applying \eqref{12} to the second term on the right-hand side of \eqref{2-5}, we have 
  \begin{align} 
    \int_{\mathbb{T}} u_h^{n + 1} (x) \, dx & = \int_{\mathbb{T}} u_h^{n} (x) \, dx + \Delta t \int_{\mathbb{T}} g^{n + 1} (x) \, dx. \label{2-3} 
  \end{align} 
  Using \eqref{2-3} recursively, we obtain 
  \begin{align} 
    \int_{\mathbb{T}} u_h^{n} (x) \, dx = \int_{\mathbb{T}} u_h^{0} (x) \, dx + \Delta t \sum_{k = 1}^{n} \int_{\mathbb{T}} g^k (x) \, dx. \label{2-7} 
  \end{align} 
  By \eqref{2-4} and $ u_h^0 = \Pi_h^{2s} u_0 $, we have 
  \begin{align} 
    \int_{\mathbb{T}} u_h^{0} (x) \, dx = \int_{\mathbb{T}} u_0 (x) \, dx. 
  \end{align} 
  Applying this to the first term on the right-hand side of \eqref{2-7}, we obtain \eqref{2-8}. 
  This completes the proof. 
\end{proof} 

\section{Proof of the main results} \label{s:3} 

\subsection{Stability} 
We discuss some properties of the operator $ \Pi_h^{2s} $. 
As a preliminary, for a positive integer $ s $ and $ v \in C^s ([0, 1]) $, we define the Hermite interpolation $ I_h^{2 s + 1} $ by 
\begin{gather} 
  (I_h^{2 s + 1} v)|_{E_j} \in \mathbb{P}_{2s + 1} (E_j) \quad \text{for} \ 1 \le j \le N_x, \\ 
  \partial_x^k (I_h^{2 s + 1} v) (x_j) = \partial_x^k v (x_j) \quad \text{for} \ 0 \le k \le s, \ 0 \le j \le N_x. 
\end{gather} 
We note that $ I_h $ is an operator for functions in $ C^s ([0, 1]) $, rather than for periodic functions. 
It is known that the operator $ I_h^{2 s + 1}: C^s ([0, 1]) \to C^s ([0, 1]) $ has the properties below: 
\begin{enumerate}[label=(I\arabic*)] 
  \item \label{i1} For any $ \varphi \in H^{s + 1} (0, 1) $, the following holds: 
  \begin{align} 
    & \int_0^1 \abs{\partial_x^{s + 1} \varphi (x)}^2 dx \\ 
    & = \int_0^1 \abs{\partial_x^{s + 1} (I_h^{2 s + 1} \varphi) (x)}^2 dx + \int_0^1 \abs{\partial_x^{s + 1} \varphi (x) - \partial_x^{s + 1} (I_h^{2 s + 1} \varphi) (x)}^2 dx. \label{3-1} 
  \end{align} 
  \item \label{i2} There exists a positive constant $ K $ depending on $ s $ such that for any $  0 \le k \le s + 1 $ 
  \begin{empheq}{alignat=2}  
    \snorm{k, 2}{\varphi - I_h^{2 s + 1} \varphi} & \leq K h^{s + 1 - k} \snorm{s + 1, 2}{\varphi} \quad & & \text{for} \ \varphi \in H^{s + 1} (0, 1), \label{3-2-2} \\ 
    \snorm{k, 2}{\varphi - I_h^{2 s + 1} \varphi} & \leq K h^{2 s + 2 - k} \snorm{2 s + 2, 2}{\varphi} \quad & & \text{for} \ \varphi \in H^{2 s + 2} (0, 1). \label{3-2} 
  \end{empheq} 
\end{enumerate} 
See Theorem~4 of \cite{67ScVa} for \ref{i1}. 
See Theorems~7 and 9 of \cite{67ScVa} for \ref{i2}. 

Here, we consider the relationship between the interpolation operator $ \Pi_h^{2s} $ and the Hermite interpolation operator. 
For $ v \in H^s (\mathbb{T}) $, let $ w $, $ V $, and $ W $ be defined by 
\begin{align}
  w = \Pi_h^{2s} v, \quad 
  V(x) = \int_0^x v(y) \, dy, \quad \text{and} \quad 
  W(x) = \int_0^x w(y) \, dy. 
\end{align}
By the definition of the interpolation operator $ \Pi_h^{2s} $, we have 
\begin{align} 
  \partial_x^k w (x_j) = \partial_x^k v (x_j) \quad \text{for} \ k  = 0, \ldots, s - 1, \ j  = 1,\ldots,N_x. \label{6-2-2} 
\end{align} 
Since $ \partial_x V = v $ and $ \partial_x W = w $, \eqref{6-2-2} implies 
\begin{align} 
  \partial_x^k V(x_j) = \partial_x^k W(x_j) \quad \text{for} \ k = 1, \ldots, s, \ j = 1, \ldots, N_x. 
  \label{6-2} 
\end{align} 
By the definition of $ \Pi_h^{2s} $, we also have 
\begin{align} 
  \int_{E_j} w(x) \, dx = \int_{E_j} v(x) \, dx \quad \text{for} \ j = 1,\ldots,N_x. 
\end{align} 
Using 
\begin{align} 
  V (x_j) & = V (0) + \sum_{i = 1}^{j} \int_{E_i} v (x)\, dx, \\ 
  W (x_j) & = W (0) + \sum_{i = 1}^{j} \int_{E_i} w (x)\, dx, 
\end{align} 
and $ V (0) = W (0) = 0 $, we have 
\begin{align} 
  V (x_j) = W (x_j) \quad \text{for} \ j = 1, \ldots, N_x. \label{6-3} 
\end{align} 
Since $ w $ is a piecewise polynomial of degree at most $ 2 s $, $ W|_{E_j} $ is a polynomial of degree at most $ 2s + 1 $. 
Combining this fact with \eqref{6-2} and \eqref{6-3}, we see that $ W $ is identical to $ I_h^{2s + 1} V $. This means that the interpolation operators $ \Pi_h^{2s} $ and $ I_h^{2 s + 1} $ satisfy 
\begin{align} 
  \Pi_h^{2s} v = \partial_x (I_h^{2s+1} V), \label{3-3} 
\end{align} 
where $ V $ is a primitive function of $ v $. 

Substituting $ \varphi = V $ into \eqref{3-1} and using $ \partial_x V = v $, we have 
\begin{align} 
  \int_0^1 \abs{\partial_x^{s} v (x)}^2 dx = \int_0^1 \abs{\partial_x^{s + 1} (I_h^{2 s + 1} V) (x)}^2 dx + \int_0^1 \abs{\partial_x^{s} v (x) - \partial_x^{s + 1} (I_h^{2 s + 1} V) (x)}^2 dx. 
\end{align} 
Using \eqref{3-3}, we have 
\begin{align} 
  \int_0^1 \abs{\partial_x^{s} v (x)}^2 dx = \int_0^1 \abs{\partial_x^{s} (\Pi_h^{2s} v) (x)}^2 dx + \int_0^1 \abs{\partial_x^{s} v (x) - \partial_x^{s} (\Pi_h^{2s} v) (x)}^2 dx. 
\end{align} 
Furthermore, substituting $ \varphi = V $ into \eqref{3-2-2} and \eqref{3-2}, and using $ \partial_x V = v $ and \eqref{3-3}, we have for $ 1 \le k \le s + 1 $ 
\begin{empheq}{alignat=2} 
  \snorm{k - 1, 2}{v - \Pi_h^{2s} v} & \leq K h^{s + 1 - k} \snorm{s, 2}{v} \quad & & \text{for} \ v \in H^s (\mathbb{T}), \\ 
  \snorm{k - 1, 2}{v - \Pi_h^{2s} v} & \leq K h^{2 s + 2 - k} \snorm{2 s + 1, 2}{v} \quad & & \text{for} \ v \in H^{2 s + 1} (\mathbb{T}). 
\end{empheq} 

From these observations, we conclude that $ \Pi_h^{2s} $ possesses the following properties: 
\begin{enumerate}[label=(P\arabic*)] 
  \item \label{P1} The following equality holds for any $ v \in H^s(\mathbb{T}) $ 
  \begin{align} 
    \snorm{s, 2}{v}^2 = \snorm{s, 2}{\Pi_h^{2s} v}^2 + \snorm{s, 2}{(I - \Pi_h^{2s}) v}^2, \label{P1-1} 
  \end{align} 
  where $ I $ is the identity operator. 
  \item \label{P2} There exists a positive constant $ K $ independent of $ h $ such that 
  \begin{empheq}{alignat=2}  
    \norm{0, 2}{v - \Pi_h^{2s} v} & \leq K h^{s} \snorm{s, 2}{v} \quad & & \text{for} \ v \in H^{s} (\mathbb{T}), \label{P2-1} \\ 
    \norm{0, 2}{v - \Pi_h^{2s} v} & \leq K h^{2 s + 1} \snorm{2 s + 1, 2}{v} \quad & & \text{for} \ v \in H^{2 s + 1} (\mathbb{T}), \label{P2-2} \\ 
    \snorm{s, 2}{v - \Pi_h^{2s} v} & \leq K h^{s + 1} \snorm{2 s + 1, 2}{v} \quad & & \text{for} \ v \in H^{2 s + 1} (\mathbb{T}). \label{P2-3} 
  \end{empheq} 
\end{enumerate} 

By \ref{P1} and \ref{P2}, we obtain the stability of the operator $ \Pi_h^{2s} $ with respect to the weighted $ H^s $-norm defined by \eqref{3-4}, as stated in the next lemma. 

\begin{lemma} \label{l:1} 
  Let $ s $ be a positive integer. Let $ \Pi_h^{2s} $ be the interpolation operator defined by \eqref{6-5} and \eqref{6-1}. Then, there exists a positive constant $ C $, depending only on $ K $, such that for any $ v \in H^s (\mathbb{T}) $, the following inequality holds: 
  \begin{align} 
    \norm{s, 2, \Delta}{\Pi_h^{2s} v} \leq (1 + C \Delta t) \norm{s, 2, \Delta}{v}. \label{l1-1} 
  \end{align} 
\end{lemma} 

\begin{proof}[Proof of Lemma~\ref{l:1}] 
  By a direct calculation, we have for any $ b_1, b_2 \in \mathbb{R} $ and $ d > 0 $ 
\begin{align} 
  (b_1 + b_2)^2 \leq (1 + d) b_1^2 + \left(1 + \frac{1}{d}\right) b_2^2. \label{6-4} 
\end{align} 
By the triangle inequality, we have 
\begin{align}
  \norm{s, 2, \Delta}{\Pi_h^{2s} v}^2 
  & = \norm{0, 2}{\Pi_h^{2s} v}^2 + \frac{h^{2s}}{\Delta t} \snorm{s, 2}{\Pi_h^{2s} v}^2 \\ 
  & \leq \left(\norm{0, 2}{v} + \norm{0, 2}{(I - \Pi_h^{2s}) v}\right)^2 + \frac{h^{2s}}{\Delta t} \snorm{s, 2}{\Pi_h^{2s} v}^2. \label{l1-2}
\end{align}
Applying \eqref{6-4} to the first term on the right-hand side of \eqref{l1-2}, we have 
\begin{align} 
  \norm{s, 2, \Delta}{\Pi_h^{2s} v}^2 
  & \leq \left(1 + K^2 \Delta t\right) \norm{0, 2}{v}^2 + \left(1 + \frac{1}{K^2 \Delta t} \right) \norm{0, 2}{(I - \Pi_h^{2s}) v}^2 + \frac{h^{2s}}{\Delta t} \snorm{s, 2}{\Pi_h^{2s} v}^2. \label{l1-3} 
\end{align} 
Substituting $ v = (I - \Pi_h^{2s}) w $ into \eqref{P2-1}, we have for $ w \in H^s (\mathbb{T}) $ 
\begin{align} 
  \norm{0, 2}{(I - \Pi_h^{2s}) w} 
  & \leq K h^s \snorm{H^s}{(I - \Pi_h^{2s}) ((I - \Pi_h^{2s}) w)} \\ 
  & = K h^s \snorm{H^s}{(I - \Pi_h^{2s}) w}, \label{l1-4} 
\end{align} 
where we used $ \Pi_h^{2s} (\Pi_h^{2s} w) = \Pi_h^{2s} w $. 
Applying \eqref{l1-4} to the second term on the right-hand side of \eqref{l1-3}, we have 
\begin{align} 
  \norm{s, 2, \Delta}{\Pi_h^{2s} v}^2 
  & \leq \left(1 + K^2 \Delta t\right) \norm{0, 2}{v}^2 + \left(1 + \frac{1}{K^2 \Delta t} \right) (K h^s)^2 \snorm{s, 2}{(I - \Pi_h^{2s}) v}^2 + \frac{h^{2s}}{\Delta t} \snorm{s, 2}{\Pi_h^{2s} v}^2 \\
  & \leq \left(1 + K^2 \Delta t\right) \norm{0, 2}{v}^2 + \left(1 + K^2 \Delta t \right)
  \frac{h^{2s}}{\Delta t} \left(
    \snorm{s, 2}{(I - \Pi_h^{2s}) v}^2 +  \snorm{s, 2}{\Pi_h^{2s} v}^2 
  \right). 
  \label{l1-5} 
\end{align} 
Applying \eqref{P1-1} to the second term on the right-hand side of \eqref{l1-5}, we have 
\begin{align}
  \norm{s, 2, \Delta}{\Pi_h^{2s} v}^2 
  & \leq \left(1 + K^2 \Delta t\right) \norm{0, 2}{v}^2 + \left(1 + K^2 \Delta t \right)
  \frac{h^{2s}}{\Delta t} \snorm{s, 2}{v}^2. 
\end{align} 
By taking $ C > 0 $ such that $ 1 + K^2 \Delta t \leq (1 + C \Delta t)^2 $, we have 
\begin{align}
  \norm{s, 2, \Delta}{\Pi_h^{2s} v}^2 
  & \leq (1 + C \Delta t)^2 
  \left(\norm{0, 2}{v}^2 + \frac{h^{2s}}{\Delta t} \snorm{s, 2}{v}^2\right) = (1 + C \Delta t)^2 \norm{s, 2, \Delta}{v}^2, 
\end{align} 
which implies \eqref{l1-1}. 
\end{proof} 

For $ \Delta t > 0 $ and $ w: \mathbb{T} \to \mathbb{R} $, we define $ X_{\Delta t, 0} [w]: \mathbb{T} \to \mathbb{T} $ by 
\begin{align} 
  X_{\Delta t, 0} [w] (x) = x - w (x) \Delta t. \label{l3-2} 
\end{align} 
In the following auxiliary lemma, we discuss the stability of the mapping $ v \mapsto v \circ X_{\Delta t, 0} [w] $ with respect to the weighted $ H^s $-norm. 

\begin{lemma} \label{l:3} 
  For a positive integer $ s $, let $ w \in W^{s, \infty} (\mathbb{T}) $. 
  Suppose that \eqref{l2-14} holds for some $ \kappa > 0 $. 
  Then there exists a positive constant $ C $, depending only on $ s $, $ \kappa $, and $ \norm{s - 1, \infty}{\partial_x w} $, such that for any $ v \in H^s (\mathbb{T}) $ and $ 0 < \Delta t \le 1 $ with $ \Delta t \snorm{1, \infty}{w} \le 1 / 2 $ 
  \begin{align} 
    \norm{s, 2, \Delta}{v \circ X_{\Delta t, 0} [w]} \leq (1 + C \Delta t) \norm{s, 2, \Delta}{v}. \label{c2-2} 
  \end{align} 
\end{lemma} 

The proof of Lemma~\ref{l:3} is given in the appendix. 

The lemma below establishes the stability of the entire scheme \eqref{3}--\eqref{4} with respect to the weighted $H^s$-norm. 

\begin{lemma} \label{l:2} 
  We assume that $ a \in C ([0, T]; W^{s + 1, \infty}(\mathbb T)) $. 
  Suppose that $ h $ and $ \Delta t $ satisfy \eqref{l2-14-1}, \eqref{l2-14}, and \eqref{l2-14-3} for $ \kappa > 0 $. 
  Then, there exists a positive constant $ C $, independent of $ h $ and $ \Delta t $, such that for any $ v \in H^s (\mathbb{T}) $, $ n = 1, \ldots, N_t $ 
  \begin{align} 
    \norm{s, 2, \Delta}{\frac{1}{2} \Pi_h^{2s} \left((v \circ X_{\Delta t, +}^{n} + v \circ X_{\Delta t, -}^{n}) Y_{\Delta t}^{n}\right)} \le (1 + C \Delta t) \norm{s, 2, \Delta}{v}. \label{l2-1} 
  \end{align} 
  Furthermore, let $ \{u_h^n\}_{n=0}^{N_t} $ be the numerical solution of the scheme \eqref{3}--\eqref{4}. 
  Then, there exists a positive constant $ C $, independent of $ h $ and $ \Delta t $, such that for $ n = 0, \ldots, N_t $ 
  \begin{align} 
    \norm{s, 2, \Delta}{u_h^n} & \le C \left(\norm{s, 2, \Delta}{u_0} + \Delta t \sum_{k = 1}^{n} \norm{s, 2, \Delta}{g^k}\right). \label{l2-2} 
  \end{align} 
\end{lemma} 

\begin{proof} 
  From the Gagliardo--Nirenberg inequality 
  \begin{align} 
    \snorm{k, 2}{v} \leq C \norm{0, 2}{v}^{1 - k / s} \snorm{s, 2}{v}^{k / s} \quad \text{for} \ 0 \le k \le s 
  \end{align} 
  and the inequality of arithmetic and geometric means 
  \begin{align} 
    b_1^{1 - w} b_2^w \leq (1 - w) b_1 + w b_2 \quad \text{for} \ b_1, b_2 > 0 \ \text{and} \ 0 \le w \le 1, 
  \end{align} 
  we deduce that there exists $ C > 0 $ depending only on $ s $ such that 
  \begin{align} 
    \snorm{k, 2}{v} \leq C (\norm{0, 2}{v} + \snorm{s, 2}{v}) \quad \text{for} \ 0 \le k \le s. \label{l2-9-1} 
  \end{align} 

  We prove that there exists a positive constant $ C $ such that the following three inequalities hold for any $ v \in H^s (\mathbb{T}) $: 
  \begin{align} 
    \norm{s, 2, \Delta}{Y_{\Delta t}^n \, v} & \leq (1 + C \Delta t) \norm{s, 2, \Delta}{v}, \label{l2-3} \\ 
    \norm{s, 2, \Delta}{v \circ X_{\Delta t, +}^n} & \leq (1 + C \Delta t) \norm{s, 2, \Delta}{v}, \label{l2-4} \\ 
    \norm{s, 2, \Delta}{v \circ X_{\Delta t, -}^n} & \leq (1 + C \Delta t) \norm{s, 2, \Delta}{v}. \label{l2-5} 
  \end{align} 

  First, we prove \eqref{l2-3}. 
  Let $ a^n (x) = a (x, t_n) $. 
  Since $ Y_{\Delta t}^n (x) = 1 - \partial_x a^n (x) \Delta t $, there exists a positive constant $ C $ depending only on $ \norm{s, \infty}{\partial_x a} $ such that 
  \begin{align} 
    \norm{0, \infty}{Y_{\Delta t}^n} & \le 1 + C \Delta t, \label{l2-6} \\ 
    \snorm{k, \infty}{Y_{\Delta t}^n} & \le C \Delta t \quad \text{for} \ 1 \le k \le s. \label{l2-7} 
  \end{align} 
  By \eqref{l2-6}, we have 
  \begin{align} 
    \norm{0, 2}{Y_{\Delta t}^n \, v} & \leq \norm{0, \infty}{Y_{\Delta t}^n} \norm{0, 2}{v} \le (1 + C \Delta t) \norm{0, 2}{v}. \label{l2-11} 
  \end{align} 
  By the triangle inequality, we have 
  \begin{align} 
    \snorm{s, 2}{Y_{\Delta t}^n \, v} & = \norm{0, 2}{\sum_{k = 0}^s \binom{s}{k} \partial_x^{s - k} Y_{\Delta t}^n \cdot \partial_x^k v} \\ 
    & \le \norm{0, \infty}{Y_{\Delta t}^n} \snorm{s, 2}{v} + \sum_{k = 0}^{s - 1} \binom{s}{k} \snorm{s - k, \infty}{Y_{\Delta t}^n} \snorm{k, 2}{v}. \label{l2-8} 
  \end{align} 
  By \eqref{l2-6} and \eqref{l2-7}, \eqref{l2-8} becomes 
  \begin{align} 
    \snorm{s, 2}{Y_{\Delta t}^n \, v} 
    & \le (1 + C \Delta t) \snorm{s, 2}{v} + C \Delta t \sum_{k = 0}^{s - 1} \snorm{k, 2}{v}. \label{l2-9-2} 
  \end{align} 
  Applying \eqref{l2-9-1} to the second term on the right-hand side of \eqref{l2-9-2}, we have 
  \begin{align}
    \snorm{s, 2}{Y_{\Delta t}^n \, v} 
    & \le (1 + C \Delta t) \snorm{s, 2}{v} + 
    C \Delta t \sum_{k = 0}^{s - 1} (\norm{0, 2}{v} + \snorm{s, 2}{v}) \\ 
    & \le (1 + C \Delta t) \snorm{s, 2}{v} + C \Delta t \norm{0, 2}{v}. \label{l2-10} 
  \end{align}
  Combining \eqref{l2-11} and \eqref{l2-10}, we obtain 
  \begin{align} 
    \norm{s, 2, \Delta}{Y_{\Delta t}^n \, v}^2 
    & \leq \left((1 + C \Delta t) \norm{0, 2}{v}\right)^2 \\ 
    & \quad + \frac{h^{2s}}{\Delta t} \left(C \Delta t \norm{0, 2}{v} + (1 + C \Delta t) \snorm{s, 2}{v}\right)^2. \label{l2-12}
  \end{align} 
  Applying \eqref{6-4} to the second term on the right-hand side of \eqref{l2-12}, we have 
  \begin{align} 
    \norm{s, 2, \Delta}{Y_{\Delta t}^n \, v}^2 & \leq (1 + C \Delta t)^2 \norm{0, 2}{v}^2 + \frac{h^{2s}}{\Delta t} (C \Delta t)^2 \left(1 + \frac{1}{\Delta t}\right) \norm{0, 2}{v}^2 \\ 
    & \quad + \frac{h^{2s}}{\Delta t} (1 + C \Delta t)^2 (1 + \Delta t) \snorm{s, 2}{v}^2. \label{l2-13} 
  \end{align} 
  By \eqref{l2-14} and $ \Delta t \le 1 $, we have 
  \begin{align} 
    \norm{s, 2, \Delta}{Y_{\Delta t}^n \, v}^2 \le (1 + C \Delta t) \left(\norm{0, 2}{v}^2 + \frac{h^{2s}}{\Delta t} \snorm{s, 2}{v}^2\right) 
    = (1 + C \Delta t) \norm{s, 2, \Delta}{v}^2, 
  \end{align} 
  which implies \eqref{l2-3}. 

  Next, we prove the inequality \eqref{l2-4}. For $ v \in H^s (\mathbb{T}) $, we define $ v_+ \in H^s (\mathbb{T}) $ by 
  \begin{align} 
    v_+ (x) = v (x + \sigma \sqrt{\Delta t}). 
  \end{align} 
  Let $ X_{\Delta t, 0} [a^n] $ be defined by 
  \begin{align}
    X_{\Delta t, 0} [a^n] (x) = x - a^n (x) \Delta t. 
  \end{align}
  Then, we have 
  \begin{align} 
    v (X_{\Delta t, +}^n (x)) = v_{+} (X_{\Delta t, 0} [a^n] (x)). 
  \end{align} 
  By Lemma~\ref{l:3}, if \eqref{l2-14-1}, \eqref{l2-14}, and \eqref{l2-14-3} are satisfied, there exists $ C $ depending only on $ s $, $ \kappa $, and $ \norm{L^\infty (0, T; W^{s - 1, \infty} (\mathbb{T}))}{\partial_x a} $ such that 
  \begin{align} 
    \norm{s, 2, \Delta}{v \circ X_{\Delta t, +}^n} & = \norm{s, 2, \Delta}{v_{+} \circ X_{\Delta t, 0} [a^n]} \\ 
    & \le (1 + C \Delta t) \norm{s, 2, \Delta}{v_{+}}. 
  \end{align} 
  Combining this with $ \norm{s, 2, \Delta}{v_+} = \norm{s, 2, \Delta}{v} $, we obtain \eqref{l2-4}. We can also obtain \eqref{l2-5} in a similar manner. 
  
  We now prove \eqref{l2-1}. 
  By Lemma~\ref{l:1}, we have 
  \begin{align} 
    & \norm{s, 2, \Delta}{\frac{1}{2} \Pi_h^{2s} \left((v \circ X_{\Delta t, +}^n + v \circ X_{\Delta t, -}^n) Y_{\Delta t}^n\right)} \\ 
    & \leq (1 + C \Delta t) \norm{s, 2, \Delta}{\frac{1}{2} (v \circ X_{\Delta t, +}^n + v \circ X_{\Delta t, -}^n) Y_{\Delta t}^n}. \label{l2-18} 
  \end{align} 
  By \eqref{l2-3}, we have 
  \begin{align}
    \norm{s, 2, \Delta}{\frac{1}{2} (v \circ X_{\Delta t, +}^n + v \circ X_{\Delta t, -}^n) Y_{\Delta t}^n}
    \leq (1 + C \Delta t) \norm{s, 2, \Delta}{\frac{1}{2} (v \circ X_{\Delta t, +}^n + v \circ X_{\Delta t, -}^n)}. \label{l2-16} 
  \end{align} 
  Furthermore, by \eqref{l2-4} and \eqref{l2-5}, we have 
  \begin{align} 
    \norm{s, 2, \Delta}{\frac{1}{2} (v \circ X_{\Delta t, +}^n + v \circ X_{\Delta t, -}^n)} \leq (1 + C \Delta t) \norm{s, 2, \Delta}{v}. \label{l2-17} 
  \end{align} 
  Combining \eqref{l2-18} and \eqref{l2-16}, we see that there exists $ C > 0 $ such that 
  \begin{align} 
    & \norm{s, 2, \Delta}{\frac{1}{2} \Pi_h^{2s} \left((v \circ X_{\Delta t, +}^n + v \circ X_{\Delta t, -}^n) Y_{\Delta t}^n\right)} \\ 
    & \le (1 + C \Delta t) \norm{s, 2, \Delta}{\frac{1}{2} (v \circ X_{\Delta t, +}^n + v \circ X_{\Delta t, -}^n)}. \label{l2-19} 
  \end{align} 
  By \eqref{l2-17} and \eqref{l2-19}, there exists $ C > 0 $ such that \eqref{l2-1} holds. 

  In the following, we prove \eqref{l2-2}. 
  Since $ \{u_h^n\}_n $ is defined by \eqref{3}--\eqref{4}, we have 
  \begin{align} 
    \norm{s, 2, \Delta}{u_h^{n + 1}} 
    & = \norm{s, 2, \Delta}{\frac{1}{2} \Pi_h^{2s} \left((u_h^n \circ X_{\Delta t, +}^{n + 1}) \, Y_{\Delta t}^{n + 1} + (u_h^n \circ X_{\Delta t, -}^{n + 1}) \, Y_{\Delta t}^{n + 1}\right) + \Pi_h^{2s} Z_{\Delta t}^{n + 1}} \\ 
    & \le \norm{s, 2, \Delta}{\frac{1}{2} \Pi_h^{2s} \left((u_h^n \circ X_{\Delta t, +}^{n + 1}) \, Y_{\Delta t}^{n + 1} + (u_h^n \circ X_{\Delta t, -}^{n + 1}) \, Y_{\Delta t}^{n + 1}\right)} \\ 
    & \quad + \norm{s, 2, \Delta}{\Pi_h^{2s} Z_{\Delta t}^{n + 1}}. \label{l2-20} 
  \end{align} 
  By \eqref{l1-1}, we estimate the second term on the right-hand side as 
  \begin{align}
    \norm{s, 2, \Delta}{\Pi_h^{2s} Z_{\Delta t}^{n + 1}} & \le (1 + C \Delta t) \norm{s, 2, \Delta}{Z_{\Delta t}^{n + 1}} \\ 
    & = (1 + C \Delta t) \Delta t \norm{s, 2, \Delta}{g^{n + 1}}. \label{l2-22}
  \end{align} 
  Applying \eqref{l2-1} to the first term on the right-hand side of \eqref{l2-20}, and applying \eqref{l2-22} to the second term, we see that there exists $ C > 0 $ such that 
  \begin{align} 
    \norm{s, 2, \Delta}{u_h^{n + 1}} 
    & \le (1 + C \Delta t) \norm{s, 2, \Delta}{u_h^n} + (1 + C \Delta t) \Delta t \norm{s, 2, \Delta}{g^{n + 1}}. \label{l2-21} 
  \end{align}
  Using \eqref{l2-21} recursively, we obtain 
  \begin{align} 
    \norm{s, 2, \Delta}{u_h^n} & \le (1 + C \Delta t)^n \norm{s, 2, \Delta}{u_h^0} + \Delta t \sum_{k = 1}^{n} (1 + C \Delta t)^{n - k + 1} \norm{s, 2, \Delta}{g^k}. \label{l2-23} 
  \end{align}
  Using \eqref{3} and \eqref{l1-1}, we have 
  \begin{align}
    \norm{s, 2, \Delta}{u_h^0} \le (1 + C \Delta t) \norm{s, 2, \Delta}{u_0}. \label{l2-24} 
  \end{align} 
  Applying \eqref{l2-24} to the first term on the right-hand side of \eqref{l2-23}, we have 
  \begin{align} 
    \norm{s, 2, \Delta}{u_h^n} & \le (1 + C \Delta t)^{n + 1} \norm{s, 2, \Delta}{u_0} + \Delta t \sum_{k = 1}^{n} (1 + C \Delta t)^{n - k + 1} \norm{s, 2, \Delta}{g^k}. \label{l2-25} 
  \end{align} 
  Since $ 1 + d \le e^d $ holds for any $ d > 0 $, we have 
  \begin{align} 
    (1 + C \Delta t)^n \le \left(1 + \frac{C T}{N_t}\right)^{N_t} 
    \le \left(e^{C T / {N_t}}\right)^{N_t} = e^{C T} \quad \text{for} \ 0 \le n \le N_t. \label{l2-26} 
  \end{align} 
  From \eqref{l2-26} and $ \Delta t \le 1 $, we obtain 
  \begin{align} 
    (1 + C \Delta t)^n \le e^{CT} (1 + C) \quad \text{for} \ 0 \le n \le N_t + 1. \label{l2-27} 
  \end{align} 
  Applying \eqref{l2-27} to \eqref{l2-25}, we obtain \eqref{l2-2}. This completes the proof. 
\end{proof} 

\subsection{Error estimates} 

Before we provide the proof of Theorem~\ref{t:1}, we estimate the truncation error of the semi-discrete scheme \eqref{14} 
\begin{align} 
  \eta^{n + 1} = u^{n + 1} - \left(\frac{1}{2} (u^n \circ X_{\Delta t, +}^{n + 1}) \, Y_{\Delta t}^{n + 1} + \frac{1}{2} (u^n \circ X_{\Delta t, -}^{n + 1}) \, Y_{\Delta t}^{n + 1} + Z_{\Delta t}^{n + 1}\right). \label{t-1} 
\end{align} 
We decompose $ \eta^{n + 1} $ as follows: 
\begin{align} 
  & \abs{\eta^{n + 1} (x)} \\ 
  & = \biggl| u^{n + 1} (x) - \frac{1}{2}\bigl(u^n(X_{\Delta t, +}^{n+1}(x)) + u^n(X_{\Delta t, -}^{n+1}(x))\bigr) Y_{\Delta t}^{n + 1} (x) - \Delta t \, g^{n + 1} (x) \biggr| \\ 
  & \le \abs{u^{n + 1} (x) - u^n (x) - \Delta t \, \partial_t u (x, t_n)} \\ 
  & \quad + \biggl|\Delta t \partial_t u (x, t_n)  + \Delta t \, \partial_x(a^{n + 1} (x) \, u^n(x)) - \frac{\sigma^2\Delta t}{2}\partial_x^2u^n(x) - \Delta t \, g^{n + 1} (x)\biggr| \\ 
  & \quad + \biggl| u^n(x) - \Delta t \, \partial_x (a^{n + 1} (x) \, u^n(x)) + \frac{\sigma^2\Delta t}{2}\partial_x^2u^n(x) \\ 
  & \quad \qquad - \frac12\bigl(u^n(X_{\Delta t, +}^{n+1}(x)) + u^n(X_{\Delta t, -}^{n + 1}(x))\bigr) Y_{\Delta t}^{n + 1} (x) \biggr|. \\ 
  & \eqqcolon \eta_1^{n + 1} (x) + \eta_2^{n + 1} (x) + \eta_3^{n + 1} (x). \label{t-11} 
\end{align} 
By Taylor's theorem, we can estimate the first term on the right-hand side of \eqref{t-11} as 
\begin{align} 
  \abs{\eta_1^{n + 1} (x)} & = \abs{\Delta t^2 \int_{0}^{1} (1 - \theta) \partial_t^2 u (x, t_{n} + \theta \Delta t) \, d\theta} \\ 
  & \le \frac{\Delta t^2}{2} \sup_{t \in (t_n, t_{n + 1})} \abs{\partial_t^2 u (x, t)}. \label{t-8} 
\end{align} 
Using \eqref{5}, we obtain the bound for $ \eta_2^{n + 1} $ as 
\begin{align} 
  & \abs{\eta_2^{n + 1} (x)} \\
  & = \biggl|\Delta t \partial_t u (x, t_n)  + \Delta t \, \partial_x (a^{n + 1} (x) \, u^n(x)) - \frac{\sigma^2\Delta t}{2}\partial_x^2u^n(x) - \Delta t \, g^{n + 1} (x) \\ 
  & \qquad - \Delta t \left(\partial_t u (x, t_n) + \partial_x (a^n (x) \, u^n (x)) - \frac{1}{2} \sigma^2 \partial_{xx} u^n (x) - g^n (x)\right)\biggr| \\ 
  & = \Delta t \abs{\left(\partial_x (a^{n + 1} (x) u^n (x)) - \partial_x (a^{n} (x) u^n (x))\right) 
  - \left(g^{n + 1} (x) - g^{n} (x)\right)} \\ 
  & = \Delta t^2 \abs{\int_{0}^{1} \partial_t \partial_x \left(a (x, t_n + \theta \Delta t) u^n (x)\right) \, d \theta - \int_{0}^{1} \partial_t g (x, t_n + \theta \Delta t) \, d \theta}, 
\end{align} 
which implies there exists $ C > 0 $ depending on 
$ \norm{W^{1, \infty} (t_n, t_{n + 1}; W^{1, \infty} (\mathbb{T}))}{a} $, 
$ \norm{L^\infty (t_n, t_{n + 1}; W^{1, \infty} (\mathbb{T}))}{u} $, and 
$ \norm{W^{1, \infty} (t_n, t_{n + 1}; L^\infty (\mathbb{T}))}{g} $ 
such that 
\begin{align} 
  \norm{0, \infty}{\eta_2^{n + 1}} \le C \Delta t^2. \label{t-6} 
\end{align} 
Let $ d_\pm = - a^{n + 1} (x) \Delta t \pm \sigma \sqrt{\Delta t} $. 
By Taylor's theorem, we have 
\begin{align} 
  & \biggl| \frac{1}{2} \left(u^n(X_{\Delta t, +}^{n+1}(x)) + u^n(X_{\Delta t, -}^{n+1}(x))\right) \\ 
  & \qquad - \left(u^n (x) - \Delta t \, a^{n + 1} (x) \partial_x u^n (x) + \frac{\sigma^2 \Delta t}{2} \partial_{xx} u^n (x) \right)\biggr| \\ 
  & = \Delta t^2 \biggl| \frac{a^{n + 1} (x)^2}{2} \partial_{xx} u^n (x) - \left(\frac{a^{n + 1} (x)^3 \Delta t}{6} + \frac{a^{n + 1} (x) \sigma^2}{2}\right) \partial_{xxx} u^n (x) \\ 
  & \quad + \frac{1}{12} \left(\frac{d_+^4}{\Delta t^2} \int_0^1 (1 - \theta)^3 \partial_x^4 u^n (x + \theta d_+) \, d\theta + \frac{d_-^4}{\Delta t^2} \int_0^1 (1 - \theta)^3 \partial_x^4 u^n (x + \theta d_-) \, d\theta\right) \biggr| \\ 
  & \le C \Delta t^2, 
\end{align} 
where $ C > 0 $ depends only on $ \norm{0, \infty}{a^{n + 1}} $, $ \norm{4, \infty}{u^n} $, and $ \sigma $. 
Therefore, we have 
\begin{align} 
  \abs{\eta_3^{n + 1} (x)} 
  & = \biggl| \frac{1}{2} \left(u^n(X_{\Delta t, +}^{n+1}(x)) + u^n(X_{\Delta t, -}^{n+1}(x))\right) Y_{\Delta t}^{n + 1} (x) \\ 
  & \qquad - \left(u^n (x) - \Delta t \, \partial_x (a^{n + 1} (x) u^n (x)) + \frac{\sigma^2 \Delta t}{2} \partial_{xx} u^n (x) \right)\biggr| \\ 
  & = \biggl| \biggl[ \frac{1}{2} \left(u^n(X_{\Delta t, +}^{n+1}(x)) + u^n(X_{\Delta t, -}^{n+1}(x))\right) \\ 
  & \qquad - \left(u^n (x) - \Delta t \, a^{n + 1} (x) \partial_x u^n (x) + \frac{\sigma^2 \Delta t}{2} \partial_{xx} u^n (x) \right) \biggr] Y_{\Delta t}^{n + 1} (x) \\ 
  & \qquad + \Delta t^2 \left(a^{n + 1} (x) \partial_x u^n (x) - \frac{\sigma^2}{2} \partial_{xx} u^n (x)\right) \partial_x a^{n + 1} (x) \biggr| \\
  & \le C \Delta t^2, \label{t-4} 
\end{align} 
where $ C > 0 $ depends only on 
$ \norm{1, \infty}{a^{n + 1}} $, 
$ \norm{4, \infty}{u^n} $, and $ \sigma $. 

By \eqref{t-8}, \eqref{t-6}, and \eqref{t-4}, we obtain 
\begin{align} 
  \norm{0, \infty}{\eta} \le 
  \norm{0, \infty}{\eta_1} + \norm{0, \infty}{\eta_2} + \norm{0, \infty}{\eta_3} 
  \le C \Delta t^2. 
\end{align} 
This implies that there exists $ C > 0 $ depending on 
$ \norm{W^{2, \infty} (0, T; W^{4, \infty} (\mathbb{T}))}{u} $, 
$ \norm{W^{1, \infty} (0, T; W^{1, \infty} (\mathbb{T}))}{a} $, 
$ \norm{W^{1, \infty} (0, T; L^\infty (\mathbb{T}))}{g} $, and 
$ \sigma $ such that 
\begin{align}
  \norm{0, 2}{\eta} \le C \Delta t^2. \label{t-9} 
\end{align} 
Furthermore, by a direct calculation, we can see that there exists $ C > 0 $ depending on 
$ \norm{W^{2, \infty} (0, T; W^{s + 4, \infty} (\mathbb{T}))}{u} $, 
$ \norm{W^{1, \infty} (0, T; W^{s + 1, \infty} (\mathbb{T}))}{a} $, 
$ \norm{W^{1, \infty} (0, T; W^{s, \infty} (\mathbb{T}))}{g} $, and 
$ \sigma $ such that 
\begin{align} 
  \snorm{s, 2}{\eta} \le C \Delta t^2. \label{t-10} 
\end{align} 
where we omit the details. 
Combining \eqref{t-9} and \eqref{t-10}, we obtain 
\begin{align} 
  \norm{s, 2, \Delta}{\eta} \leq C \left(1 + \frac{h^s}{\sqrt{\Delta t}}\right) \Delta t^2. \label{t1} 
\end{align} 

In the rest of this section, we give the proof of Theorem~\ref{t:1}. 

\begin{proof}[Proof of Theorem~\ref{t:1}] 
  Let $ u^n (x) = u(x, t_n) $ and $ e_h^n = u^n - u_h^n $. 
  First, we estimate $ e_h^0 $. 
  By \eqref{P2-2}, \eqref{P2-3}, and $ e_h^0 = u_0 - \Pi_h^{2s} u_0 $, we have 
  \begin{align} 
    \norm{0, 2}{e_h^0} & \le K h^{2 s + 1} \snorm{2 s + 1, 2}{u_0}, \\ 
    \snorm{s, 2}{e_h^0} & \le K h^{s + 1} \snorm{2 s + 1, 2}{u_0}. 
  \end{align} 
  By these two inequalities, we have 
  \begin{align} 
    \norm{s, 2, \Delta}{e_h^0}^2 & \le 
    \norm{0, 2}{e_h^0}^2 + \frac{h^{2 s}}{\Delta t} \snorm{s, 2}{e_h^0}^2 \\ 
     & \le \left(1 + \frac{1}{\Delta t}\right) K^2 h^{4 s + 2} \snorm{2 s + 1, 2}{u_0}^2. 
  \end{align} 
  Since $ \Delta t \le 1 $, there exists $ C > 0 $ such that 
  \begin{align}
    \norm{s, 2, \Delta}{e_h^0}^2 & \le C \frac{h^{4 s + 2}}{\Delta t} \snorm{2 s + 1, 2}{u_0}^2. \label{1-18} 
  \end{align}

  Next, suppose that $ 0 \le n \le N_t - 1 $. 
  By \eqref{4}, we have 
  \begin{align} 
    e_h^{n + 1} & = u^{n + 1} - u_h^{n + 1} \\ 
    & = u^{n + 1} - \Pi_h^{2s} \biggl( \frac{1}{2} (u_h^n \circ X_{\Delta t, +}^{n + 1}) \, Y_{\Delta t}^{n + 1} + \frac{1}{2} (u_h^n \circ X_{\Delta t, -}^{n + 1}) \, Y_{\Delta t}^{n + 1} + Z_{\Delta t}^{n + 1} \biggr). 
  \end{align}
  Using $ u_h^n = u^n - e_h^n $, we have 
    \begin{align}
    e_h^{n + 1} 
    & = \Pi_h^{2s}
    \left(\frac{1}{2} \left(e_h^n \circ X_{\Delta t, +}^{n + 1} + e_h^n \circ X_{\Delta t, -}^{n + 1}\right) Y_{\Delta t}^{n + 1}\right) \\ 
    & \quad + \Pi_h^{2s} \biggl( u^{n + 1} - \Bigl( \frac{1}{2} (u^n \circ X_{\Delta t, +}^{n + 1}) \, Y_{\Delta t}^{n + 1} + \frac{1}{2} (u^n \circ X_{\Delta t, -}^{n + 1}) \, Y_{\Delta t}^{n + 1} + Z_{\Delta t}^{n + 1}\Bigr) \biggr) \\ 
    & \quad  + (I - \Pi_h^{2s}) u^{n + 1}. \label{1-1-2} 
  \end{align} 
  Let $ \tilde{e}_h^n = \frac{1}{2} (e_h^n \circ X_{\Delta t, +}^{n + 1} + e_h^n \circ X_{\Delta t, -}^{n + 1}) Y_{\Delta t}^{n + 1} $ and $ \rho_h^n = u^n - \Pi_h^{2s} u^{n} $. 
  Then, \eqref{1-1-2} can be rewritten as 
  \begin{align} 
    e_h^{n + 1} & = \Pi_h^{2s} \tilde{e}_h^n + \Pi_h^{2s} \eta^{n + 1} + (I - \Pi_h^{2s}) u^{n + 1}, \label{1-1} 
  \end{align} 
  where we recall that $ \eta^n $ is defined by \eqref{t-1}. 
  Taking $ L^2 $-norm of both sides of \eqref{1-1}, we have 
  \begin{align} 
    \norm{0, 2}{e_h^{n + 1}}^2 
    & = \norm{0, 2}{\Pi_h^{2s} \tilde{e}_h^n + \Pi_h^{2s} \eta^{n + 1} + \rho_h^{n + 1}}^2 \\ 
    & \leq \left(\norm{0, 2}{\Pi_h^{2s} \tilde{e}_h^n} + \norm{0, 2}{\Pi_h^{2s} \eta^{n + 1}} + \norm{0, 2}{\rho_h^{n + 1}}\right)^2. \label{1-2} 
  \end{align} 
  Applying the inequality 
  \begin{align} 
    (a_1 + a_2 + a_3)^2 \leq (1 + 2 \Delta t) a_1^2 + \left(2 + \frac{1}{\Delta t}\right) (a_2^2 + a_3^2) \quad \text{for} \ a_1, a_2, a_3 \in \mathbb{R}, \ \Delta t > 0 
  \end{align} 
  to \eqref{1-2}, we have 
  \begin{align} 
    \norm{0, 2}{e_h^{n + 1}}^2 & \leq (1 + 2 \Delta t) \norm{0, 2}{\Pi_h^{2s} \tilde{e}_h^n}^2 + \left(2 + \frac{1}{\Delta t}\right) \left(\norm{0, 2}{\Pi_h^{2s} \eta^{n + 1}}^2 + \norm{0, 2}{\rho_h^{n + 1}}^2\right). \label{1-7} 
  \end{align} 
  Applying $ \Pi_h^{2s} $ to both sides of \eqref{1-1}, we have 
  \begin{align} 
    \Pi_h^{2s} e_h^{n + 1} & = \Pi_h^{2s} \tilde{e}_h^n + \Pi_h^{2s} \eta^{n + 1}. \label{1-3} 
  \end{align} 
  Applying $ I - \Pi_h^{2s} $ to both side of \eqref{1-1}, we have 
  \begin{align} 
    (I - \Pi_h^{2s}) e_h^{n + 1} = (I - \Pi_h^{2s}) u^{n + 1}. \label{1-4} 
  \end{align} 
  Substituting $ v = e_h^{n + 1} $ into \eqref{P1-1}, and using \eqref{1-3} and \eqref{1-4}, we obtain 
  \begin{align}
    \snorm{s, 2}{e_h^{n + 1}}^2 & = \snorm{s, 2}{\Pi_h^{2s} \tilde{e}_h^n + \Pi_h^{2s} \eta^{n + 1}}^2 + \snorm{s, 2}{\rho_h^{n + 1}}^2 \\ 
    & \leq \left(\snorm{s, 2}{\Pi_h^{2s} \tilde{e}_h^n} + \snorm{s, 2}{\Pi_h^{2s} \eta^{n + 1}}\right)^2 + \snorm{s, 2}{\rho_h^{n + 1}}^2. \label{1-5} 
  \end{align} 
  Applying the inequality 
  \begin{align} 
    (a_1 + a_2)^2 \leq (1 + \Delta t) a_1^2 + \left(1 + \frac{1}{\Delta t}\right)a_2^2 \quad \text{for} \ a_1, a_2 \in \mathbb{R} 
  \end{align} 
  to the first term of the right-hand side of \eqref{1-5}, we obtain 
  \begin{align}
    \snorm{s, 2}{e_h^{n + 1}}^2 & \leq (1 + \Delta t) \snorm{s, 2}{\Pi_h^{2s} \tilde{e}_h^n}^2 + \left(1 + \frac{1}{\Delta t}\right) \snorm{s, 2}{\Pi_h^{2s} \eta^{n + 1}}^2 + \snorm{s, 2}{\rho_h^{n + 1}}^2. \label{1-6} 
  \end{align} 
  Applying \eqref{1-7} and \eqref{1-6} to $ \norm{s, 2, \Delta}{e_h^{n + 1}}^2 = \norm{0, 2}{e_h^{n + 1}}^2 + \frac{h^{2s}}{\Delta t} \snorm{s, 2}{e_h^{n + 1}}^2 $, we obtain 
  \begin{align} 
    \norm{s, 2, \Delta}{e_h^{n + 1}}^2 & \leq 
    (1 + 2 \Delta t) \norm{0, 2}{\Pi_h^{2s} \tilde{e}_h^n}^2 + 
    \left(2 + \frac{1}{\Delta t}\right) \left(\norm{0, 2}{\Pi_h^{2s} \eta^{n + 1}}^2 + \norm{0, 2}{\rho_h^{n + 1}}^2\right) \\ 
    & \quad + (1 + \Delta t) \frac{h^{2s}}{\Delta t} \snorm{s, 2}{\Pi_h^{2s} \tilde{e}_h^n}^2 + \left(1 + \frac{1}{\Delta t}\right) \frac{h^{2s}}{\Delta t} \snorm{s, 2}{\Pi_h^{2s} \eta^{n + 1}}^2 + \frac{h^{2s}}{\Delta t} \snorm{s, 2}{\rho_h^{n + 1}}^2 \\ 
    & \leq (1 + 2 \Delta t) \norm{s, 2, \Delta}{\Pi_h^{2s} \tilde{e}_h^n}^2 + \left(2 + \frac{1}{\Delta t}\right) \norm{s, 2, \Delta}{\Pi_h^{2s} \eta^{n + 1}}^2 \\ 
    & \quad + \left(2 + \frac{1}{\Delta t}\right) \norm{0, 2}{\rho_h^{n + 1}}^2
    + \frac{h^{2s}}{\Delta t} \snorm{s, 2}{\rho_h^{n + 1}}^2. \label{1-9} 
  \end{align} 
  We estimate the second, third, and fourth term on the right-hand side of \eqref{1-9}. 
  Substituting $ v = \eta^{n + 1} $ into \eqref{l1-1}, we have 
  \begin{align} 
    \norm{s, 2, \Delta}{\Pi_h^{2s} \eta^{n + 1}} & \leq (1 + C \Delta t) \norm{s, 2, \Delta}{\eta^{n + 1}}. \label{1-10} 
  \end{align} 
  Combining \eqref{l2-14} and \eqref{t1}, we have 
  \begin{align} 
    \norm{s, 2, \Delta}{\eta^{n + 1}} \le C \Delta t^2. \label{1-10-2} 
  \end{align} 
  Combining \eqref{1-10-2}, \eqref{1-10}, and $ \Delta t \le 1 $, we see there exists $ C > 0 $ such that 
  \begin{align} 
    \norm{s, 2, \Delta}{\Pi_h^{2s} \eta^{n + 1}} & \leq C \Delta t^2. \label{1-8} 
  \end{align} 
  Furthermore, by \eqref{P2-2} and \eqref{P2-3}, we have 
  \begin{align} 
    \norm{0, 2}{\rho_h^{n + 1}} & = \norm{0, 2}{(I - \Pi_h^{2s}) u^{n + 1}} \leq K h^{2s + 1} \norm{2 s + 1, 2}{u^{n + 1}}, \label{1-11} \\ 
    \snorm{s, 2}{\rho_h^{n + 1}} & = \snorm{s, 2}{(I - \Pi_h^{2s}) u^{n + 1}} \leq K h^{s + 1} \norm{2 s + 1, 2}{u^{n + 1}}. \label{1-12} 
  \end{align} 
  Applying \eqref{1-8}, \eqref{1-11}, and \eqref{1-12} to \eqref{1-9}, we obtain 
  \begin{align} 
    \norm{s, 2, \Delta}{e_h^{n + 1}}^2 & \le
    (1 + 2 \Delta t) \norm{s, 2, \Delta}{\Pi_h^{2s} \tilde{e}_h^n}^2 + C \left(\Delta t^3 + \frac{h^{4 s + 2}}{\Delta t}\right). \label{1-14} 
  \end{align} 
  
  By Lemma~\ref{l:2}, we have 
  \begin{align}
    \norm{s, 2, \Delta}{\Pi_h^{2s} \tilde{e}_h^n} 
    & = \norm{s, 2, \Delta}{\frac{1}{2} \Pi_h^{2s} \left( (e_h^{n} \circ X_{\Delta t, +}^{n + 1} + e_h^{n} \circ X_{\Delta t, -}^{n + 1}) Y_{\Delta t}^{n + 1}\right)} \\ 
    & \le (1 + C \Delta t) \norm{s, 2, \Delta}{e_h^{n}}. \label{1-15-2} 
  \end{align} 
  Applying \eqref{1-15-2} to the first term of \eqref{1-14}, we have 
  \begin{align} 
    \norm{s, 2, \Delta}{e_h^{n + 1}}^2 & \le (1 + 2 \Delta t) (1 + C \Delta t)^2 \norm{s, 2, \Delta}{e_h^{n}}^2 + C \left(\Delta t^3 + \frac{h^{4 s + 2}}{\Delta t}\right) \\ 
    & \le (1 + C \Delta t) \norm{s, 2, \Delta}{e_h^{n}}^2 + C \left(\Delta t^3 + \frac{h^{4 s + 2}}{\Delta t}\right),  \label{1-16} 
  \end{align} 
  where we used $ \Delta t \le 1 $. 
  Using \eqref{1-16} recursively, we have for $ n \ge 1 $ 
  \begin{align} 
    \norm{s, 2, \Delta}{e_h^{n}}^2 
    & \le (1 + C \Delta t)^n \norm{s, 2, \Delta}{e_h^0}^2 + C \sum_{k = 0}^{n - 1} \left(1 + C \Delta t\right)^k \left(\Delta t^3 + \frac{h^{4 s + 2}}{\Delta t}\right). \label{1-17} 
  \end{align} 
  Using \eqref{l2-26}, we have 
  \begin{align} 
    \norm{s, 2, \Delta}{e_h^{n}}^2 & \le 
    e^{C T} \norm{s, 2, \Delta}{e_h^0}^2 + C e^{C T} \sum_{k = 0}^{n - 1} \left(\Delta t^3 + \frac{h^{4 s + 2}}{\Delta t}\right). \label{1-19} 
  \end{align} 
  Applying \eqref{1-18} to the first term on the right-hand side of \eqref{1-19}, we have 
  \begin{align} 
    \norm{s, 2, \Delta}{e_h^{n}}^2 & \le C e^{C T} \frac{h^{4 s + 2}}{\Delta t} \snorm{2 s + 1, 2}{u_0}^2 + C n e^{C T} \left(\Delta t^3 + \frac{h^{4 s + 2}}{\Delta t}\right). 
  \end{align} 
  Using $ n \le \frac{T}{\Delta t} $ and $ \Delta t \le 1 $, we have 
  \begin{align}
    \norm{s, 2, \Delta}{e_h^{n}}^2 & 
    \le C \left(\frac{h^{4 s + 2}}{\Delta t} + \Delta t^2 + \frac{h^{4 s + 2}}{\Delta t^2}\right) \\ 
    & \le C \left(\Delta t^2 + \frac{h^{4 s + 2}}{\Delta t^2}\right), \label{1-22-2} 
  \end{align}
  which implies the weighted $ H^s $-error estimate 
  \begin{align} 
    \norm{s, 2, \Delta}{e_h^{n}} \le C \left(\Delta t + \frac{h^{2 s + 1}}{\Delta t}\right). \label{1-22} 
  \end{align} 

  By \eqref{3-4}, we have for $ v \in H^s (\mathbb{T}) $ 
  \begin{align} 
    \norm{0, 2}{v} & \le \norm{s, 2, \Delta}{v}, \label{1-20} \\ 
    \snorm{s, 2}{v} & \le \frac{\sqrt{\Delta t}}{h^s} \norm{s, 2, \Delta}{v}. \label{1-21} 
  \end{align} 
  From \eqref{1-22} and \eqref{1-20}, we deduce \eqref{t1-1}. 
  Furthermore, \eqref{l2-9-1} implies that 
  \begin{align} 
    \norm{s, 2}{v} \le C (\norm{0, 2}{v} + \snorm{s, 2}{v}). \label{1-23} 
  \end{align} 
  By \eqref{1-20}, \eqref{1-21}, and \eqref{1-23}, there exists $ C > 0 $ depending only on $ s $ such that for any $ v \in H^s (\mathbb{T}) $ 
  \begin{align} 
    \norm{s, 2}{v} \le C \left(1 + \frac{\sqrt{\Delta t}}{h^s}\right) \norm{s, 2, \Delta}{v}. 
  \end{align} 
  By \eqref{l2-14}, there exists $ C > 0 $ depending only on $ s $ and $ \kappa $ such that 
  \begin{align} 
    \norm{s, 2}{v} \le C \frac{\sqrt{\Delta t}}{h^s} \norm{s, 2, \Delta}{v}. \label{1-24} 
  \end{align} 
  From \eqref{1-22} and \eqref{1-24}, we deduce \eqref{t1-2}. 
  This completes the proof. 
\end{proof} 

\section{Numerical Experiments} \label{s:e} 
\begin{figure}[htbp] 
\centering 
\begin{tikzpicture} 
\begin{axis}[ 
  width=0.5\textwidth, height=0.5\textwidth,
  xmode=log, ymode=log,
  xlabel={$h$}, ylabel={relative $ L^2 $-error}, 
  xtick={0.1,0.05,0.025,0.0125},
  xticklabels={$\frac{1}{10}$,$\frac{1}{20}$,$\frac{1}{40}$,$\frac{1}{80}$},
  ytick={1e-1,1e-2,1e-3,1e-4,1e-5,1e-6,1e-7,1e-8},
  mark size=4pt,
  legend style={font=\scriptsize, at={(1.03,1)}, anchor=north west, draw=black, fill=white}, legend cell align={left},
  grid=both, grid style={gray!50}, minor grid style={gray!15}, tick style={gray!50}, tick label style={font=\scriptsize}, label style={font=\scriptsize}, title style={font=\scriptsize},
]
\addplot[black, mark=triangle, solid] coordinates {
  (0.1,2.1317e-01) (0.05,1.7108e-01) (0.025,1.1420e-01) (0.0125,6.6657e-02) 
};
\addlegendentry{(i) FSL degree 1}
\addplot[black, mark=o, solid] coordinates {
  (0.1,1.4444e-04) (0.05,1.4118e-05) (0.025,2.4282e-06) (0.0125,5.3270e-07) 
};
\addlegendentry{(ii) FSL degree 3} 
\addplot[black, mark=x, solid] coordinates {
  (0.1,2.8102e-04) (0.05,2.4152e-05) (0.025,3.1578e-06) (0.0125,5.8164e-07) 
};
\addlegendentry{(iii) Conservative FSL degree 2}
\addplot[black, mark=square, solid] coordinates {
  (0.1,9.1031e-06) (0.05,3.7130e-07) (0.025,3.6046e-08) (0.0125,3.4795e-08) 
};
\addlegendentry{(iv) Conservative FSL degree 4}
\end{axis} 
\end{tikzpicture} 
\caption{Log-log plot of the relative $L^2$-error versus spatial mesh size $h$, at fixed time step $\Delta t = 10^{-5}$. (i) the FSL scheme with linear interpolation \eqref{7-2}--\eqref{8-2};  (ii) the FSL scheme with cubic Hermite interpolation;  (iii) the conservative FSL scheme with the interpolation $ \Pi_h^2 $ (\eqref{3}--\eqref{4} with $ s = 1 $);  (iv) the conservative FSL scheme with the interpolation $ \Pi_h^4 $ (\eqref{3}--\eqref{4} with $ s = 2 $) } 
\label{5-7} 
\end{figure}
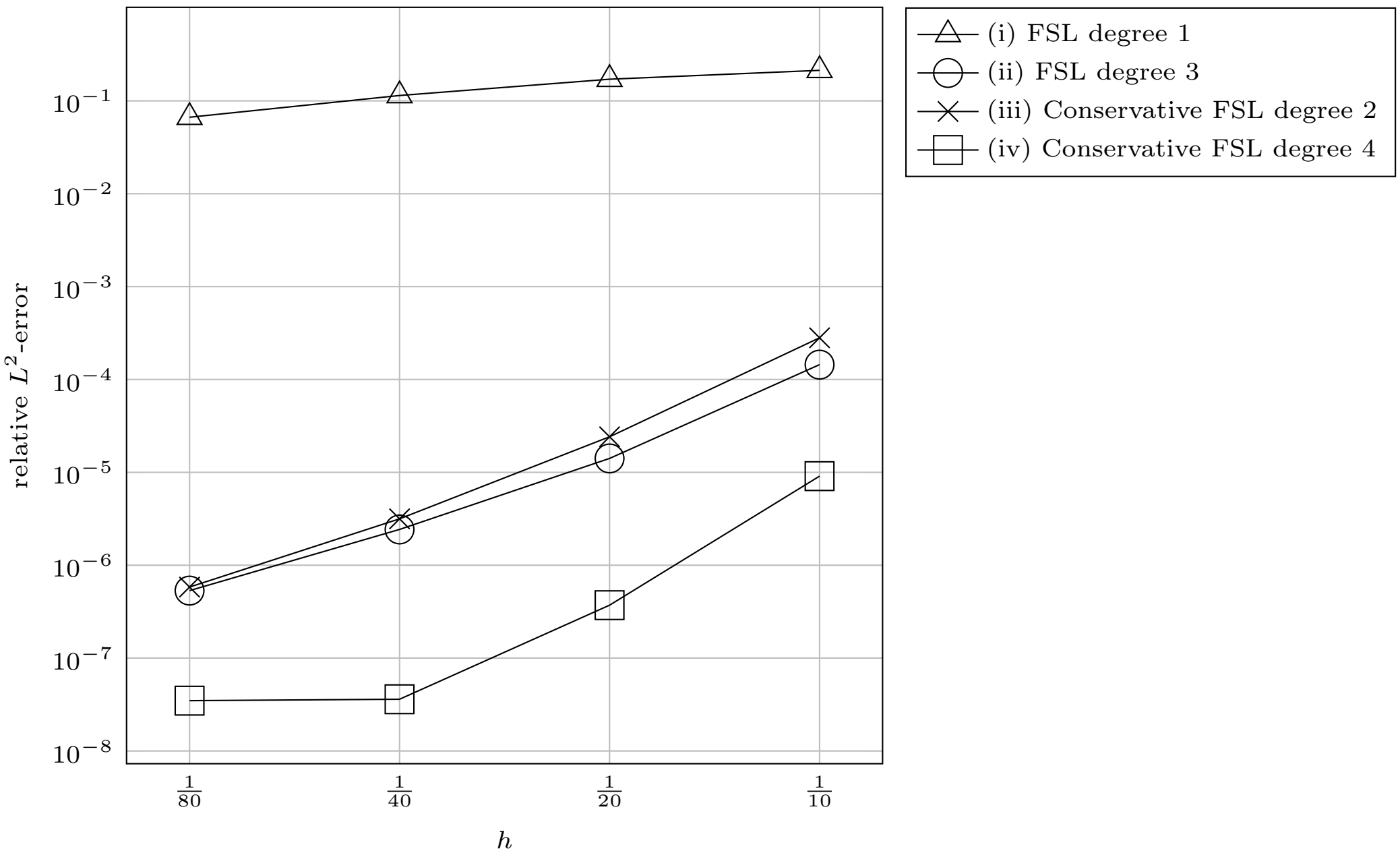

We conduct numerical tests with the following settings: 
\begin{align} 
  \frac{1}{2}\sigma^2 & = \frac{1}{1000}, \label{5-2} \\ 
  a (x, t) & = 1 + \frac{1}{10} \sin (2\pi x), \label{5-3} \\
  g (x, t) & = \frac{3\pi}{10} \cos (2\pi x) + \frac{\pi}{10} \, e^{- \frac{1}{250} \pi^2 t} \sin \big(2\pi (2x - t) \big), \label{5-5} \\ 
  u_0 (x) & = \frac{3}{2} + \frac{1}{2} \sin (2 \pi x). \label{5-6} 
\end{align} 
Then, the exact solution to \eqref{5}--\eqref{6} can be written as 
\begin{align} 
  u (x, t) = \frac{3}{2} + \frac{1}{2} \, e^{- \frac{1}{250} \pi^2 t} \sin \big(2\pi (x - t) \big). \label{5-4} 
\end{align} 
We compare four schemes: 
\begin{enumerate}[label=(\roman*)] 
  \item The scheme \eqref{7-2}--\eqref{8-2}, namely, the FSL scheme with linear interpolation. 
  \item The FSL scheme with cubic Hermite interpolation obtained by replacing $ I_h^{\mathrm{lin}} $ in \eqref{7-2}--\eqref{8-2} by the cubic Hermite interpolation operator. 
  \item The scheme \eqref{3}--\eqref{4} with $ s = 1 $, namely, the conservative FSL scheme with the interpolation operator $ \Pi_h^2 $. 
  \item The scheme \eqref{3}--\eqref{4} with $ s = 2 $, namely, the conservative FSL scheme with the interpolation operator $ \Pi_h^4 $. 
\end{enumerate} 
Let the relative $ L^2 $-error be defined by $ \| u^{N_t} - u_h^{N_t} \|_{0, 2} / \| u^{N_t} \|_{0, 2} $. 
We fix $ \Delta t = 10^{-5} $, $ T = 0.1 $, and set $ N_x = 10, 20, 40, 80 $. 
The results are shown in Figure~\ref{5-7}. 

While the scheme (iv) achieves the highest accuracy as expected, a theoretical explanation for why the schemes (ii) and (iii) exhibit nearly identical accuracy, despite the different degrees of their interpolation operators, remains to be found. 

\begin{figure}[htbp]
\centering
\begin{tikzpicture}
\begin{axis}[
  width=0.5\textwidth, height=0.5\textwidth,
  xmode=linear, ymode=log,
  xtick={0,5,10,15,20}, 
  xlabel={Time $ t_n $}, ylabel={Relative mass error $ e_{\mathrm{mass}}^n $}, 
  ymin=1e-16, ymax=1e0,
  ytick={1e0,1e-1,1e-2,1e-3,1e-4,1e-5,1e-6,1e-7,1e-8,1e-9,1e-10,1e-11,1e-12,1e-13,1e-14,1e-15,1e-16},
  yticklabels={$10^{0}$,,,,$10^{-4}$,,,,$10^{-8}$,,,,$10^{-12}$,,,,$10^{-16}$},
  mark size=4pt,
  legend style={font=\scriptsize, at={(1.03,1)}, anchor=north west, draw=black, fill=white}, legend cell align={left},
  grid=both, grid style={gray!50}, minor grid style={gray!15}, tick style={gray!50}, tick label style={font=\scriptsize}, label style={font=\scriptsize}, title style={font=\scriptsize}, 
] 
\addplot[black, mark=triangle, solid] coordinates {
  (0.0,3.7547e-03) (2.0,3.8397e-03) (4.0,4.1277e-03) (6.0,4.4549e-03) (8.0,4.7896e-03) (10.0,5.1256e-03) (12.0,5.4618e-03) (14.0,5.7978e-03) (16.0,6.1338e-03) (18.0,6.4697e-03) (20.0,6.8054e-03) 
};
\addlegendentry{(i) FSL degree 1}
\addplot[black, mark=o, solid] coordinates {
  (0.0,3.7547e-03) (2.0,3.6648e-03) (4.0,3.6581e-03) (6.0,3.6514e-03) (8.0,3.6447e-03) (10.0,3.6381e-03) (12.0,3.6315e-03) (14.0,3.6250e-03) (16.0,3.6184e-03) (18.0,3.6119e-03) (20.0,3.6054e-03) 
};
\addlegendentry{(ii) FSL degree 3}
\addplot[black, mark=x, solid] coordinates {
  (0.0,2.7705e-16) (2.0,1.1096e-13) (4.0,2.2247e-13) (6.0,3.3315e-13) (8.0,4.4536e-13) (10.0,5.5659e-13) (12.0,6.6907e-13) (14.0,7.8059e-13) (16.0,8.9168e-13) (18.0,1.0036e-12) (20.0,1.1154e-12) 
};
\addlegendentry{(iii) Conservative FSL degree 2}
\addplot[black, mark=square, solid] coordinates {
  (0.0,2.7705e-16) (2.0,1.2038e-13) (4.0,2.3286e-13) (6.0,3.5379e-13) (8.0,4.6711e-13) (10.0,6.0051e-13) (12.0,6.8196e-13) (14.0,8.0109e-13) (16.0,9.2673e-13) (18.0,1.0377e-12) (20.0,1.1480e-12) 
};
\addlegendentry{(iv) Conservative FSL degree 4}
\end{axis}
\end{tikzpicture} 
\caption{Relative mass error $ e_{\mathrm{mass}}^n $ defined by \eqref{5-12} versus time $ t_n $ for the schemes (i)--(iv) } 
\label{5-11} 
\end{figure}
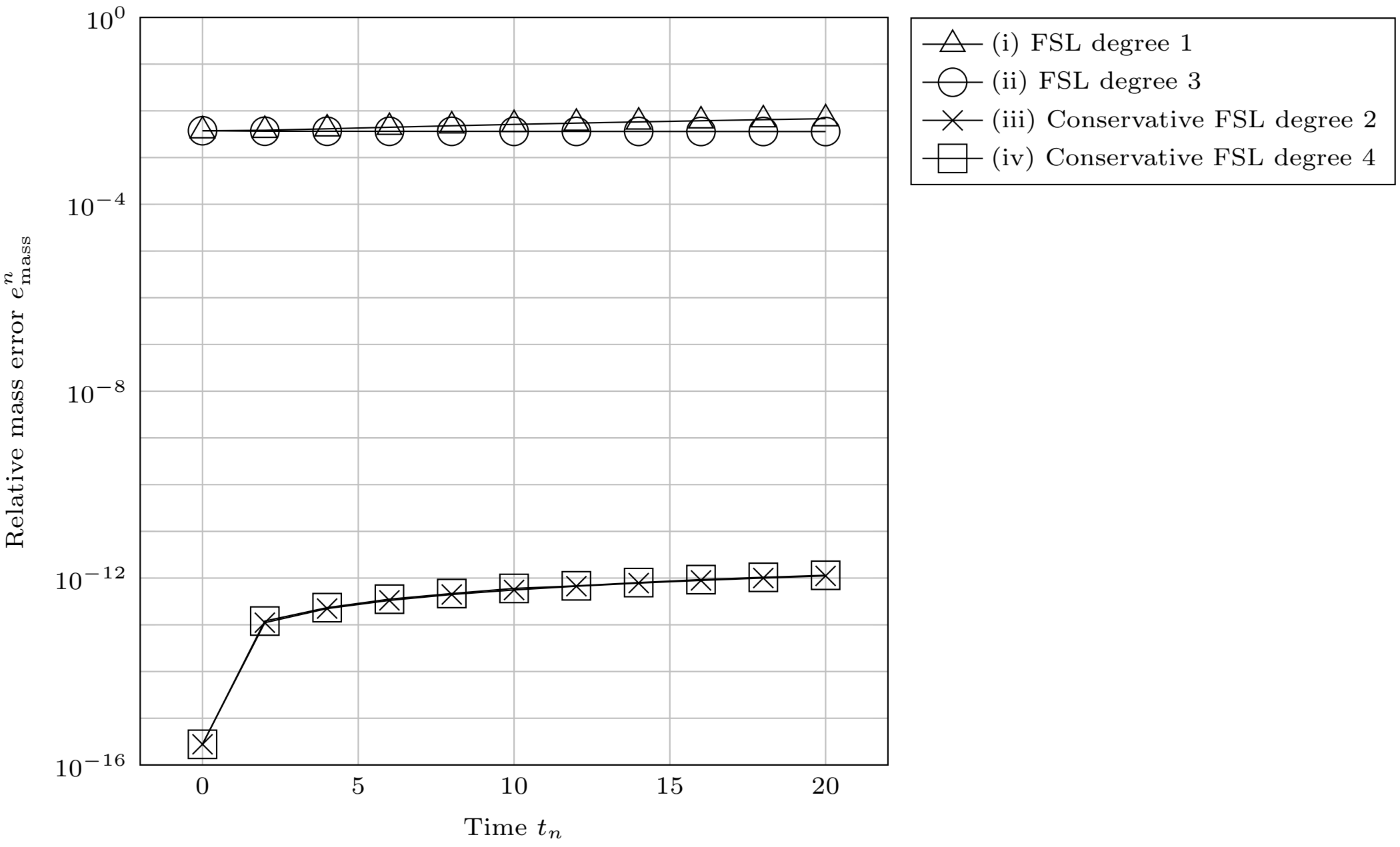

Next, we verify numerically the mass-conservation property by solving the problem \eqref{5}--\eqref{6} with the settings \eqref{5-2}, \eqref{5-3}, $ u_0 = \exp (\sin (2\pi x) + \cos (20 \pi x)) $ and $ g \equiv 0 $. We use the uniform mesh in space with $ h = 1 / 32 $, and set time step size $ \Delta t = 10^{-3} $. 
We define the relative mass error at $ t_n $ by 
\begin{align} 
  e_{\mathrm{mass}}^n = \frac{\abs{\int_\mathbb{T} u^n (x) \, dx - \int_\mathbb{T} u_h^n (x) \, dx}}{\int_\mathbb{T} u^n (x) \, dx}. \label{5-12}
\end{align} 
The results are shown in Figure~\ref{5-11}. 
While the schemes (i) and (ii) yield mass errors $e_{\mathrm{mass}}^{N_t}$ larger than $10^{-3}$, the two conservative FSL schemes achieve significantly smaller errors below $10^{-11}$. 

\appendix 
\section{Auxiliary lemmas} \label{s:a} 

We recall Fa\`a di Bruno's formula (see \cite{02Jo}) 
\begin{align} 
  \partial_x^s (v (w (x))) & = \sum_{k = 0}^s  \sum_{\alpha \in \mathcal{A}_{s,k}} C_{s,\alpha} \prod_{j = 1}^s \left(\partial_x^j w(x)\right)^{\alpha_j} \partial_x^k v(w(x)), \label{l3-5} 
\end{align} 
where 
\begin{gather} 
  \mathcal{A}_{s, k} = \biggl\{\alpha = (\alpha_1, \ldots, \alpha_s) \in (\mathbb{N} \cup \{0\})^s \biggm| \sum_{j=1}^{s} j \alpha_j = s, \ \sum_{j=1}^{s} \alpha_j = k\biggr\}, \label{6-7} \\ 
  C_{s,\alpha} = \frac{s!}{\prod_{j=1}^{s}\alpha_j! \, (j!)^{\alpha_j}}. \label{6-6}
\end{gather} 

\begin{proof}[Proof of Lemma~\ref{l:3}] 
  We prove the following two inequalities: 
  \begin{align} 
    \norm{0, \infty}{\partial_y ({X_{\Delta t, 0}} [w])^{-1}}^{1 / 2} & \le 1 + \norm{0, \infty}{\partial_x w} \Delta t, \label{l3-4} \\ 
    \norm{0, 2}{v \circ X_{\Delta t, 0} [w]} & \le (1 + \norm{0, \infty}{\partial_x w} \Delta t) \norm{0, 2}{v} \quad \text{for} \ v \in L^2 (\mathbb{T}). \label{c8} 
  \end{align} 

  Differentiating \eqref{l3-2}, we have 
  \begin{align} 
    \partial_x X_{\Delta t, 0} [w] (x) = 1 - \partial_x w (x) \Delta t. \label{l3-3} 
  \end{align} 
  Combining $ \Delta t \snorm{1, \infty}{w} \leq 1 / 2 $ with this, we have 
  \begin{align} 
    \frac{1}{2} \leq \partial_x X_{\Delta t, 0} [w] \leq \frac{3}{2}. 
  \end{align} 
  We denote the inverse of $ X_{\Delta t, 0} [w] $ by $ ({X_{\Delta t, 0}} [w])^{-1} $. 
  These bounds imply that $ X_{\Delta t,0}[w] $ is a bi-Lipschitz function. 
  Therefore, by \eqref{l3-3}, we have for almost every $ y \in \mathbb{T} $ 
  \begin{align} 
    \partial_y ({X_{\Delta t, 0}} [w])^{-1} (y) & = \frac{1}{1 - \partial_x w ((X_{\Delta t, 0} [w])^{-1} (y)) \Delta t}. 
  \end{align}
  Since the inequality 
  \begin{align} 
    (1 - a)^{- 1 / 2} \leq 1 + a 
  \end{align} 
  holds for any $ a \in [0, 1/2] $, we have 
  \begin{align} 
    \left(\partial_y ({X_{\Delta t, 0}} [w])^{-1} (y)\right)^{1 / 2} 
    & \le \left(\frac{1}{1 - \abs{\partial_x w ((X_{\Delta t, 0} [w])^{-1} (y))} \Delta t}\right)^{1 / 2} \\ 
    & \le 1 + \abs{\partial_x w ((X_{\Delta t, 0} [w])^{-1} (y))} \Delta t, 
  \end{align} 
  which implies \eqref{l3-4}. 

  Next, we prove \eqref{c8}. Using the change of variables $ y = X_{\Delta t, 0} [w] (x) $, we have 
  \begin{align} 
    \norm{0, 2}{v \circ {X_{\Delta t, 0}} [w]}^2 
    & = \int_{\mathbb{T}} v ({X_{\Delta t, 0}} [w] (x))^2 \, dx \\ 
    & = \int_{\mathbb{T}} v (y)^2 | \partial_y ({X_{\Delta t, 0}} [w])^{-1} (y) | \, dy \\ 
    & \leq \norm{0, \infty}{\partial_y ({X_{\Delta t, 0}} [w])^{-1}} \norm{0, 2}{v}^2, \label{a12-2} 
  \end{align} 
  which implies 
  \begin{align}
    \norm{0, 2}{v \circ {X_{\Delta t, 0}} [w]} \leq \norm{0, \infty}{\partial_y ({X_{\Delta t, 0}} [w])^{-1}}^{1 / 2} \norm{0, 2}{v}. \label{m4} 
  \end{align} 
  Applying \eqref{l3-4} to \eqref{m4}, we obtain \eqref{c8}. 

  In the following, we estimate the seminorm $ \snorm{s, 2}{v \circ X_{\Delta t, 0} [w]} $. 
  By \eqref{l3-5},  we have 
  \begin{align} 
    \snorm{s, 2}{v \circ X_{\Delta t, 0} [w]} 
    & = \norm{0, 2}{\partial_x^s (v \circ X_{\Delta t, 0} [w])} \\ 
    & = \norm{0, 2}{\sum_{k = 0}^s  \sum_{\alpha \in \mathcal{A}_{s,k}} C_{s,\alpha} \prod_{j = 1}^s \left(\partial_x^j X_{\Delta t, 0} [w]\right)^{\alpha_j} \cdot (\partial_x^k v) \circ X_{\Delta t, 0} [w]} \\ 
    & \leq
    \sum_{\alpha \in \mathcal{A}_{s,s}} C_{s,\alpha} \prod_{j = 1}^s \snorm{j - 1, \infty}{\partial_x X_{\Delta t, 0} [w]}^{\alpha_j}\norm{0, 2}{(\partial_x^s v) \circ X_{\Delta t, 0} [w]} \\ 
    & \quad + \sum_{k = 0}^{s - 1} \sum_{\alpha \in \mathcal{A}_{s,k}} C_{s,\alpha} \prod_{j = 1}^s \snorm{j - 1, \infty}{\partial_x X_{\Delta t, 0} [w]}^{\alpha_j}\norm{0, 2}{(\partial_x^k v) \circ X_{\Delta t, 0} [w]}. \label{l3-6} 
  \end{align} 
  By \eqref{l3-3}, we have 
  \begin{align} 
    \norm{0, \infty}{\partial_x X_{\Delta t, 0} [w]} & \le 1 + \Delta t \norm{0, \infty}{\partial_x w}, \label{l3-7} \\ 
    \snorm{k, \infty}{\partial_x X_{\Delta t, 0} [w]} & = \Delta t \snorm{k, \infty}{\partial_x w} \quad \text{for} \ 1 \le k \le s - 1. \label{l3-8} 
  \end{align} 
  Furthermore, \eqref{l3-8} implies 
  \begin{align} 
    \snorm{k, \infty}{\partial_x X_{\Delta t, 0} [w]} & \le \Delta t \norm{s - 1, \infty}{\partial_x w} \quad \text{for} \ 1 \le k \le s - 1. \label{l3-8-2} 
  \end{align} 
  Substituting $ k = s $ into \eqref{6-7}, we see that $ \mathcal{A}_{s, s} = \{(s, 0,\ldots, 0)\} $. By \eqref{6-6}, we also have $ C_{s, (s, 0,\ldots, 0)} = 1 $. Therefore, the first term on the right-hand side of \eqref{l3-6} can be bounded as 
  \begin{align} 
    & \sum_{\alpha \in \mathcal{A}_{s, s}} C_{s,\alpha} \prod_{j = 1}^s \snorm{j - 1, \infty}{\partial_x X_{\Delta t, 0} [w]}^{\alpha_j}\norm{0, 2}{(\partial_x^s v) \circ X_{\Delta t, 0} [w]} \\ 
    & = \norm{0, \infty}{\partial_x X_{\Delta t, 0} [w]}^s \norm{0, 2}{(\partial_x^s v) \circ X_{\Delta t, 0} [w]} \\ 
    & \le (1 + \Delta t \norm{0, \infty}{\partial_x w})^s \norm{0, 2}{(\partial_x^s v) \circ X_{\Delta t, 0} [w]}. 
  \end{align} 
  Using \eqref{c8}, we obtain 
  \begin{align} 
    & \sum_{\alpha \in \mathcal{A}_{s, s}} C_{s,\alpha} \prod_{j = 1}^s \snorm{j - 1, \infty}{\partial_x X_{\Delta t, 0} [w]}^{\alpha_j}\norm{0, 2}{(\partial_x^s v) \circ X_{\Delta t, 0} [w]} \\ 
    & \le (1 + \Delta t \norm{0, \infty}{\partial_x w})^{s + 1} \norm{0, 2}{\partial_x^s v} \\ 
    & \le (1 + C \Delta t) \snorm{s, 2}{v}, \label{l3-11} 
  \end{align} 
  where $ C $ depends only on $ s $ and $ \norm{0, \infty}{\partial_x w} $. 
  
  We consider the second term on the right-hand side of \eqref{l3-6}. 
  By \eqref{l3-7} and \eqref{l3-8-2}, we have 
  \begin{align} 
    \prod_{j = 1}^s \snorm{j - 1, \infty}{\partial_x X_{\Delta t, 0} [w]}^{\alpha_j} 
    & = \norm{0, \infty}{\partial_x X_{\Delta t, 0} [w]}^{\alpha_1} \prod_{j = 2}^s \snorm{j - 1, \infty}{\partial_x X_{\Delta t, 0} [w]}^{\alpha_j} \\ 
    & \le (1 + \Delta t \norm{0, \infty}{\partial_x w})^{\alpha_1} \prod_{j = 2}^s (\Delta t \norm{s - 1, \infty}{\partial_x w})^{\alpha_j} \\ 
    & \le (1 + \Delta t \norm{0, \infty}{\partial_x w})^{\alpha_1}  (\Delta t \norm{s - 1, \infty}{\partial_x w})^{\sum_{j = 2}^s \alpha_j}. \label{l3-9} 
  \end{align} 
  Assume that $ k < s $ and $ (\alpha_1, \ldots, \alpha_s) \in \mathcal{A}_{s,k} $. 
  Then, there exists $ 2 \le j \le s $ such that $ \alpha_j \ge 1 $, which implies that $ \sum_{j = 2}^s \alpha_j \ge 1 $. 
  By \eqref{l3-9} and this inequality, there exists $ C > 0 $ depending only on $ s $ and $ \norm{s - 1, \infty}{\partial_x w} $ such that 
  \begin{align} 
    \prod_{j = 1}^s \snorm{j - 1, \infty}{\partial_x X_{\Delta t, 0} [w]}^{\alpha_j} \le C \Delta t. \label{s3-10} 
  \end{align} 
  By \eqref{c8} and \eqref{s3-10}, we have for $ 0 \le k \le s - 1 $ 
  \begin{align} 
    \prod_{j = 1}^s \snorm{j - 1, \infty}{\partial_x X_{\Delta t, 0} [w]}^{\alpha_j}\norm{0, 2}{(\partial_x^k v) \circ X_{\Delta t, 0} [w]} & \le C \Delta t \norm{0, 2}{\partial_x^k v}. 
  \end{align} 
  Furthermore, by \eqref{l2-9-1}, we have 
  \begin{align}
    \prod_{j = 1}^s \snorm{j - 1, \infty}{\partial_x X_{\Delta t, 0} [w]}^{\alpha_j}\norm{0, 2}{(\partial_x^k v) \circ X_{\Delta t, 0} [w]} & \le C \Delta t (\norm{0, 2}{v} + \snorm{s, 2}{v}). \label{l3-12} 
  \end{align} 
  Applying \eqref{l3-11} and \eqref{l3-12} to the first and second terms on the right-hand side of \eqref{l3-6}, respectively, we obtain 
  \begin{align}
    \snorm{s, 2}{v \circ X_{\Delta t, 0} [w]} 
    & \le (1 + C \Delta t) \snorm{s, 2}{v} + C \Delta t (\norm{0, 2}{v} + \snorm{s, 2}{v}) \\ 
    & \le (1 + C \Delta t) \snorm{s, 2}{v} + C \Delta t \norm{0, 2}{v}. 
  \end{align} 
  Substituting $ b_1 = (1 + C \Delta t) \snorm{s, 2}{v} $, $ b_2 = C \Delta t \norm{0, 2}{v} $, and $ d = \Delta t $ into \eqref{6-4}, we obtain 
  \begin{align}
    \snorm{s, 2}{v \circ X_{\Delta t, 0} [w]}^2 
    & \le \left((1 + C \Delta t) \snorm{s, 2}{v} + C \Delta t \norm{0, 2}{v}\right)^2 \\ 
    & \le (1 + \Delta t) \left((1 + C \Delta t) \snorm{s, 2}{v}\right)^2 + \left(1 + \frac{1}{\Delta t}\right) (C \Delta t \norm{0, 2}{v})^2, 
  \end{align} 
  which implies 
  \begin{align}
    \snorm{s, 2}{v \circ X_{\Delta t, 0} [w]}^2 
    & \le (1 + C \Delta t) \snorm{s, 2}{v}^2 + C \Delta t \norm{0, 2}{v}^2. \label{l3-13}
  \end{align} 
  By \eqref{c8} and \eqref{l3-13}, we have 
  \begin{align}
    \norm{s, 2, \Delta}{v \circ X_{\Delta t, 0} [w]}^2 
    & = \norm{0, 2}{v \circ X_{\Delta t, 0} [w]}^2 + \frac{h^{2 s}}{\Delta t} \snorm{s, 2}{v \circ X_{\Delta t, 0} [w]}^2 \\ 
    & \leq (1 + C \Delta t) \norm{0, 2}{v}^2 + \frac{h^{2 s}}{\Delta t} \left((1 + C \Delta t) \snorm{s, 2}{v}^2 + C \Delta t \norm{0, 2}{v}^2\right).  
  \end{align} 
  By \eqref{l2-14}, there exists $ C > 0 $, depending only on $ s $, $ \kappa $, and $ \norm{s - 1, \infty}{\partial_x w} $, such that 
  \begin{align} 
    \norm{s, 2, \Delta}{v \circ X_{\Delta t, 0} [w]}^2 
    & \leq (1 + C (1 + \kappa) \Delta t) \norm{0, 2}{v}^2 + (1 + C \Delta t) \frac{h^{2 s}}{\Delta t} \snorm{s, 2}{v}^2 \\ 
    & \leq (1 + C \Delta t) \left(\norm{0, 2}{v}^2 + \frac{h^{2 s}}{\Delta t} \snorm{s, 2}{v}^2\right) \\ 
    & = (1 + C \Delta t) \norm{s, 2, \Delta}{v}^2, 
  \end{align} 
  which concludes the proof. 
\end{proof} 

\section*{Statements and Declarations} 
\subsection*{Funding} 
This work was supported by the Japan Society for the Promotion of Science KAKENHI Grant Number JP24KJ0964. 

\providecommand{\bysame}{\leavevmode\hbox to3em{\hrulefill}\thinspace}
\providecommand{\MR}{\relax\ifhmode\unskip\space\fi MR }
\providecommand{\MRhref}[2]{%
  \href{http://www.ams.org/mathscinet-getitem?mr=#1}{#2}
}
\providecommand{\href}[2]{#2}

\end{document}